\documentclass[
final
]{dmtcs-episciences}

\usepackage[utf8]{inputenc}
\usepackage[english]{babel}
\usepackage{amsmath}
\usepackage{amsthm}
\usepackage{amssymb}
\usepackage{amsfonts}
\usepackage{mathtools}
\usepackage{dsfont}
\usepackage{graphicx}
\usepackage{floatflt,epsfig}
\usepackage{lineno}
\usepackage{cleveref}
\usepackage{subcaption}

\usepackage{array,tabularx,tabulary,booktabs}
\usepackage{arydshln}
\usepackage{tabularray}

\newcommand{\currentrowstyle}{}
\newcolumntype{^}{>{\currentrowstyle}}

\newtheorem{lemma}{Lemma}

\newtheorem{theorem}{Theorem}

\newtheorem{proposition}{Proposition}
\newtheorem{statement}{Statement}
\newtheorem{problem}{Problem}

\newtheorem{conjecture}{Conjecture}
\newtheorem{example}{Example}

\DeclareMathOperator{\Sym}{\mathrm{Sym}}
\DeclareMathOperator{\Alt}{\mathrm{Alt}_n}

\usepackage[round]{natbib}

\author[Blanco, Golubyatnikov, Konstantinova, Maslova, and Nikiforov]{
Sa\'ul A.~Blanco\affiliationmark{1}\and Mikhail P.~Golubyatnikov\affiliationmark{2,3}\\\and Elena V.~Konstantinova\affiliationmark{4,5,6}\thanks{Supported by the Mathematical Center in Akademgorodok under the agreement No. 075-15-2025-348 with the Ministry of Science and Higher Education of the Russian Federation}\and Natalia V.~Maslova\affiliationmark{2,3}\\\and Luka A.~Nikiforov\affiliationmark{6}
}
\title{Generating the symmetric group by three prefix reversals}
\affiliation{
  Department of Computer Science, Indiana University, Bloomington, USA\\
  Krasovskii Institute of Mathematics and Mechanics Yekaterinburg, Russia\\
  Ural Mathematical Center, Yekaterinburg, Russia\\
  Sobolev Institute of Mathematics, Novosibirsk, Russia\\
  Novosibirsk State University, Novosibirsk, Russia\\
  Three Gorges Mathematical Research Center, China Three Gorges University,  Yichang,  China}
\begin{document}
\publicationdata{vol. 28:1, Permutation Patterns 2025}{2026}{7}{10.46298/dmtcs.16975}{2025-11-24; 2025-11-24; 2026-09-03}{2026-09-13}
\maketitle

\begin{abstract}
The cubic pancake graphs are Cayley graphs over the symmetric group $\mathrm{Sym}_n$ generated by three prefix reversals. There is the following open problem: characterize all the sets of three prefix reversals that generate $\mathrm{Sym}_n$. As the largest prefix reversal of length $n$ is always included in a triple, we give a complete solution of the problem when any of the two smallest or the two largest lengths but $n$ are included in a triple of prefix reversals. Moreover, some conditions implying a triple of prefix reversals does not generate $\mathrm{Sym}_n$ are considered. Computational results on the diameter and the girth of some cubic pancake graphs are presented, and conjectures for future research are formulated.
\end{abstract}


\section{Introduction}\label{sec0}

A classical and celebrated result~\cite{D05} states that two random permutations of a set of $n$ elements almost surely generate either the symmetric or the alternating group of degree $n$, and similar results hold when more generators are considered. In particular, the symmetric group can always be generated by the transposition $(1\ 2)$ together with the $n$-cycle $(1\ 2\ 3 \ldots n)$. Two well-known generating sets of $n-1$ transpositions are the Coxeter generators $\{(i\ (i+1)) \mid 1 \leqslant i < n\}$ and the star generators $\{(1\ i) \mid 1<i\leqslant n\}$ which give rise to two Cayley graphs known respectively as the \emph{bubble-sort graph} and the \emph{star graph} (see~\cite{J63,F00,AK89,K08}). In convex geometry, the bubble-sort graph is known as a $1$-skeleton of the permutahedron~\cite{T06}. The Coxeter generators play an important role in computer science~\cite{J63,K}, knot theory~\cite{BB05}, and the theory of reflection groups and Coxeter groups~\cite{H92}. The above two generating sets of transpositions were considered in~\cite{KL75,S78} where more general result was proved. Namely, it was shown that if an arbitrary  generating set of the symmetric group of degree $n$ is given by transpositions then one can generate all $n!$ permutations with no repetitions which leads to a hamiltonian cycle in the corresponding Cayley graph.

The symmetric group $\mathrm{Sym}_n$ is also generated by a set of $n-1$ involutions known as \emph{prefix reversals}, defined by
\begin{equation} \label{e1}
r_i=\begin{pmatrix} 
1 & 2 & \ldots & i-1 & i & i+1 & \ldots & n \\ 
i & i-1 & \ldots & 2 & 1 & i+1 & \ldots & n
\end{pmatrix},
\end{equation}
where $1<i\leqslant n$. The action of $r_i$ is to reverse the order of the first $i$ elements while leaving the remaining $n-i$ elements fixed. Equivalently, $r_i$ is the product of $\lfloor i/2 \rfloor$ disjoint transpositions, $r_i=(1\ i)(2\ (i-1))\cdots(\lfloor i/2\rfloor\ (i+1-\lfloor i/2\rfloor))$.

The Cayley graph generated by the set $\{r_2,r_3,\ldots,r_n\}$ is called the \emph{pancake graph}, denoted by $\mathcal{P}_n$. The pancake graph is famous because of the long-standing open problem, known as the \emph{pancake problem} which asks for the diameter of $\mathcal{P}_n$~\cite{DW75}. The \emph{girth} of a graph $\Gamma$ refers to the length of the shortest cycle of $\Gamma$, and the \emph{circumference} of $\Gamma$ refers to the length of the largest cycle in $\Gamma$. Moreover, $\Gamma$ is said to be \emph{weakly pancyclic} if $\Gamma$ has a cycle of length $\ell$ for all $g(\Gamma)\leqslant\ell\leqslant c(\Gamma)$, where $g(\Gamma)$ and $c(\Gamma)$ refer to the girth and circumference of $\Gamma$, respectively. Since the pancake graph is weakly pancyclic~\cite{KF95,STC06} (in fact, a larger family of graphs that include $\mathcal{P}_n$, called the \emph{generalized pancake graphs}, are known to be weakly pancyclic~\cite{BBP19, BB23}), this problem remains difficult to resolve. 

Our interest lies in describing sets of three prefix reversals that generate the entire symmetric group $\mathrm{Sym}_n$. 

The importance of fixed-degree pancake graphs, in particular, cubic pancake graphs as models of networks was shown in~\cite{BS03}. The authors of that work considered cubic pancake graphs as spanning subgraphs of the pancake graph $\mathcal{P}_n$ and they identified the following combinatorial necessary conditions for a triple of distinct prefix reversals $\{r_i,r_j,r_k\}$ to generate the symmetric group $\Sym_n$:  
\begin{enumerate}
    \item[$(1)$] One of the prefix reversals must be $r_n$, since otherwise no element can be moved to or from position $n$. 
    \item[$(2)$] One of the prefix reversals, say $r_i$, must satisfy $\lfloor n/2 \rfloor < i < n$, since otherwise if we take $g\in \langle r_n, r_j, r_i\rangle$ with $i<j \leqslant n/2$, and if $u <v \leqslant j$ then either $g(u) \leqslant n/2$ and $g(v) \leqslant n/2$ or $g(u) > n/2$ and $g(v)> n/2$. Thus, in this case, $2$-transitivity is lost, while $\Sym_n$ is known to be $2$-transitive in its natural action for each $n$. Moreover, if $n$ is odd, then $i>\lceil n/2\rceil$ since otherwise $\lceil n/2\rceil$ will remain fixed. 
    \item[$(3)$] At least one prefix reversal must have an even index, since otherwise no element can be moved from an odd-numbered position to an even-numbered position. 
    \item[$(4)$] At least one prefix reversal must be an odd permutation, since otherwise the generated group consists only of even permutations. 
    \item[$(5)$] The indices must satisfy the condition $\gcd(i,j,k)=1$. Otherwise, if $\gcd(i,j,k)=\ell>1$ then the elements of $\{1,2,\ldots,n\}$ are partitioned into blocks of $\ell$ consecutive positions that can never be separated. In particular, at least one index should be odd.
\end{enumerate}
\medskip
\noindent {\bf Remark. {\it Originally, in~{\rm\cite{BS03}} the authors formulate Observation $(5)$ as follows{\rm:} The generating set must contain at least two prefix reversals with relatively prime indices. Of course, if this holds then $\gcd(i,j,k)=1$. The converse does not hold. Indeed, if $i=15$, $j=12$ and $k=10$ then $\gcd(i,j)=3$, $\gcd(j,k)=2$, and $\gcd(i,k)=5$, however $\gcd(i,j,k)=1$. Moreover, it is easy to check that $\langle r_{15}, r_{12}, r_{10}\rangle = \Sym_{15}$.}
}
\medskip

The above necessary conditions are not sufficient, however, they can be used to help guide the search for generators. In particular, in~\cite{BS03} it was shown that the symmetric group is generated by the following six sets: $\{r_n,r_{n-1},r_2\}$ for any $n \geqslant 4$; $\{r_n,r_{n-1},r_3\}$ for any odd $n \geqslant 5$; $\{r_n,r_{n-1},r_{n-2}\}$ for any $n \geqslant 4$; $\{r_n,r_{n-1},r_{n-3}\}$ for any odd $n \geqslant 4$; $\{r_n,r_{n-2},r_3\}$ for any even $n \geqslant 6$; and $\{r_n,r_{n-2},r_{n-3}\}$ for any $n \geqslant 5$. The authors have used the following approach. For each of the sets above, it was shown that the set can simulate the three generators of the shuffle-exchange permutation network generated by the following three permutations written in cycle notation: $( 1\ 2\ \ldots \ n-1\ n)$ (right shuffle), $(n\ n-1\ \ldots \ 2\ 1)$ (left shuffle), and $(1\ 2)$ (transposition of the first two elements). Using the same approach, it was shown in~\cite{Gun22} that the set $\{r_n,r_{n-2},r_2\}$ for any odd $n \geqslant 5$ is also a generating set of the symmetric group.  

The set $\{r_n,r_{n-1},r_{n-2}\}$ is known as Big-3 flips, and the corresponding cubic pancake graph generated by this set is called the Big-3 pancake network~\cite{SW16} for which it was conjectured that there exists a cyclic Gray code. The girths and the diameters of the cubic pancake graphs generated by the above seven sets were studied in~\cite{KS22} and~\cite{CayleyPy}.

\subsection{Our contributions}

\medskip

In this paper, we focus on the following three problems.

\begin{problem}\label{Prob1}
Are there other generating sets of three prefix reversals for the symmetric group?
\end{problem}

\begin{problem}\label{Prob2}
What are other necessary conditions for the set of three prefix reversals to generate the symmetric group?
\end{problem}

\begin{problem}\label{Prob3}
What are sufficient conditions for the set of three prefix reversals to generate the symmetric group?
\end{problem}

In addressing Problem~\ref{Prob2}, we prove the following assertions.\\

\noindent {\bf Proposition 1.} {\it  Let $2 \leqslant k < n-l < n$ and $n>3k+1$. If $l \geqslant 2k+1$ then $ \langle r_n, r_{n-l}, r_k\rangle < \Sym_n$.}

\medskip

\noindent {\bf Proposition 2.} {Let $n \geqslant 6$ and $l \geqslant k+1$. If  $l$ divides $n+1$ then   $\langle r_n, r_{n-l},r_k\rangle < \Sym_n$.}

\medskip

In addressing the above problems, we obtain the following four classification results.

\medskip

\noindent {\bf Theorem~\ref{thm:2}.}  {\it The set $\{r_n,r_m,r_2\}$, with $2<m<n$ and $n\geqslant 4$, generates $\Sym_n$ if and only if one of the following conditions holds{\rm:}
\begin{enumerate}[$(i)$]
    \item $n$ is even and $m=n-1${\rm;} 
    \item $n$ is odd and 
    \begin{enumerate}[$(a)$]
        \item $m\in\{n-1,n-2\}$ if $n\equiv2\pmod 3$ or if $n\equiv0\pmod 3${\rm;}
        \item $m\in\{n-3,n-2,n-1\}$ if $n\equiv 1\pmod 3$. 
    \end{enumerate}
\end{enumerate}
}

\medskip

\noindent {\bf Theorem~\ref{thm:n-1}.} {\it The set $\{r_n, r_{n-1}, r_k\}$, with $2\leqslant k\leqslant n-2$, generates $\Sym_n$ if and only if one of the following conditions holds:
	\begin{enumerate}[$(i)$]
		\item If $n$ is even and $k$ is even;
		\item $n \equiv 3 \pmod{4}$ and $2\leqslant k\leqslant n-2$;
		\item $n \equiv 1 \pmod{4}$ and $k \equiv 2, 3 \pmod{4}$. 
	\end{enumerate}
}

\medskip

\noindent {\bf Theorem~\ref{thm:3}.} {\it    The set $\{r_n,r_m,r_3\}$, with $2\leqslant m<n$, $m \not =3$ and $n\geqslant 5$, generates $\Sym_n$ if and only if one of the following conditions holds{\rm:}
    \begin{enumerate}[$(i)$]
        \item $n$ is even and $m=n-2${\rm;}
        \item $m\in\{n-3,n-1\}$ and $n\equiv 1,5\pmod 6${\rm;}
        \item $m=n-1$ if $n\equiv3\pmod 6$.
    \end{enumerate}
}

\medskip    

\noindent {\bf Theorem~\ref{thm:n-2}.} {\it  The set $\{ r_n, r_{n-2}, r_k\}$, with $2\leqslant k<n-2$ and $n\geqslant5$, generates $\Sym_n$ if and only if $k$ and $n$ are of different parity.}

 \medskip

\textbf{Narrative of the results.}  The four classification theorems above provide a complete answer when one of the indices is close to either end of the range. For a generating triple containing $r_2$ or $r_3$, the other non-$r_n$ prefix reversal must have one of the few lengths closest to $n$: only $r_{n-1}$, $r_{n-2}$, and, in certain congruence classes, $r_{n-3}$ can occur. The exact list is determined by $n$  modulo $6$. If the triple contains $r_{n-1}$ then generation is determined entirely by congruence conditions on $n$ and the third index modulo $4$. If the triple contains $r_{n-2}$ the criterion is particularly simple: the third index and $n$ must have opposite parity. Thus, in each of these four broad families, the question of whether three prefix reversals generate $\mathrm{Sym}_n$ reduces to an explicit arithmetic test. The complete conditions are stated formally in Theorems~1 through 4.

\medskip

\textbf{Organization of the paper.} In Section~\ref{sec:prelim}, we describe group theory techniques, including some developed by the authors, that provide necessary and sufficient conditions for some triple $\{r_n,r_m,r_k\}$ to generate $\Sym_n$, in particular, we prove Propositions~1 and~2. In Section~\ref{sec:results}, we present our main results relating to triples of prefix reversals that generate $\Sym_n$, in particular, we prove Theorems~1 through 4.   
Moreover, we formulate conjectures in Section~\ref{sec:comp} for future research. More specifically, Conjecture~1 concerns the asymptotic number of generating triples, Conjectures~2 and~3 describe the apparent shift phenomena, Conjecture~4 describes the triangular non-generating region, and Conjecture~5 concerns the number of possible middle indices when the smallest reversal length is fixed. We finish with Section~\ref{sec:graph} where we discuss computational results relating to the cubic pancake graphs such as their diameters, girths, and hamiltonicity.

\section{Preliminaries}\label{sec:prelim}

The following notation is used. Let $\Omega$ denote the set $\{1,2,\ldots,n\}$, and let $\Sym_n$ and $\Alt$ denote the symmetric group and the alternating group of degree $n$, respectively. We use two notations for a permutation $\pi$, namely, a cyclic notation for an $n$-cycle $\pi=(\pi_1\pi_2\ldots\pi_n)$, where $\pi(1) \rightarrow  \pi(2)$, $\pi(2) \rightarrow  \pi(3)$, $\ldots$ , $\pi(n) \rightarrow  \pi(1)$, and a one-line notation $\pi=[\pi_1\pi_2\ldots\pi_n]$, where $\pi(i)=\pi_i$ for any $1\leqslant i\leqslant n$. For $\sigma,\pi\in\Sym_n$, we compose them using right-to-left multiplication, so $\sigma\pi(i)=\sigma(\pi(i))$ for any  $1\leqslant i\leqslant n$. It is easy to see that by~(\ref{e1}) we have $r_{n}r_{n-1}=(1\ 2 \ldots n-1\ n)$. Further, if it necessary in our proofs we consider cyclic permutations as a product of appropriate prefix reversals.

We now proceed to describe some  tools from group theory that are used in the proofs of our results. Let $G$ be a finite group acting on a finite set $\Omega$. The \textit{$G$-orbit of $x\in\Omega$}, or simply the \textit{orbit of $x$}, denoted by $\mathcal{O}(x)$, refers to the set $\mathcal{O}(x)=\{g(x) : g \in G\}$. A \textit{$G$-invariant set} is a union of $G$-orbits. The action is said to be {\it transitive} on $\Omega$ if for each $x, y \in \Omega$ there is $g\in G$ such that $g(x)=y$. A nonempty subset $B\subseteq \Omega$ is called a \emph{block of imprimitivity} (or simply a \emph{block}) for the action if for every $g\in G$ either $g(B) = B$ or $g(B)\cap B=\varnothing$, where $g(B)=\{g(b):b\in B\}$. Moreover, the family $\{g(B) : g\in G\}$ of translates of a block $B$ forms a $G$-invariant partition of $\Omega$, called a \emph{block system}. There are at least two block systems for every action: $\{x\}$ for $x\in \Omega$  and $\Omega$ itself, and one refers to these as \emph{trivial blocks}. The action of $G$ on $\Omega$ is \emph{primitive} if it is transitive and the only blocks are the trivial ones. Otherwise, a transitive action is said to be \emph{imprimitive}. Additionally, if $H$ is a proper subgroup of $G$, we write $H<G$.

{\bf Main ideas in the proofs.} The proofs use two recurring strategies. 

1) To show that a triple does not generate $\mathrm{Sym}_n$, we identify an invariant set of positions or a nontrivial block system, thereby proving that the generated group is intransitive or imprimitive while $\mathrm{Sym}_n$ is known to be primitive in the natural action on the set $\Omega$.

2) To prove that a triple does generate $\mathrm{Sym}_n$, we construct explicit cycles, transpositions, or $3$-cycles as words in the prefix reversals and then apply standard generation criteria for the symmetric and alternating groups. 

This division between obstruction arguments and explicit generation
arguments provides a common structure for the case analysis in Section~\ref{sec:results}.

In our proofs, we utilize a few results from group theory regarding generating sets of $\Sym_n$.

\begin{lemma}[{\rm\cite[Theorem~13.3]{W64}}]\label{WieLemma} Let $H \leqslant \text{Sym}(\Omega)=\Sym_n$ be a permutation group acting naturally on the set $\Omega=\{1, \ldots, n\}$ such that there exists a transposition $h \in H$. Then $H=\text{Sym}(\Omega)$ if and only if $H$ is primitive on $\Omega$. 
\end{lemma}

The following fact is well-known, see, for example,~\cite[Theorem~1.1]{IradmusaTaleb}.

\begin{lemma}\label{lem:Gen_(2,n)} $H=\langle (1\ 2 \ldots \ n),(a\ b)\rangle=\Sym_n$, $1 \leqslant a < b \leqslant n$, if and only if $\gcd(b - a, n) = 1$.
\end{lemma}

A more general result relating to generating sets of $\Sym_n$ formed by one $n$-cycle and one transposition is proved in~\cite{Janusz}. 

\begin{lemma}[{\rm\cite[Corollary 4]{Janusz}}]\label{lem:janusz} Let $\sigma$ be an $n$-cycle and $\tau=(a \ b)$ be a transposition. Let $q$ be an integer satisfying $\sigma^q(a)=b$. Then $\langle \sigma,\tau\rangle=\Sym_n$ if and only if $\gcd(q,n)=1$.
\end{lemma}

The following result proved in~\cite{IradmusaTaleb} is also used in our proofs. 

\begin{lemma}[{\rm\cite[Theorem~1.3]{IradmusaTaleb}}]\label{lem:Gen_(3,n)}
    Let $G$ be the permutation subgroup of $\Sym_n$ generated by the $n$-cycle $ (1\ 2\ \ldots \ n)$ and the $3$-cycle $(a\ b\ c)$ with $1\leqslant a < c < b \leqslant n$. Then, 
    
    $(i)$ If $n$ is even, then $G=\Sym_n$ if and only if $\gcd(b-a,c-a,n)=1$. 
    
    $(ii)$ If $n$ is odd, then 
    $G=\Alt$ if and only if $\gcd(b-a,c-a,n)=1$. 
\end{lemma}

We furthermore develop tools to detect when the action of $\langle r_n,r_m,r_k\rangle$ on $\Omega$ is intransitive, from which it follows that $\langle r_n,r_m,r_k\rangle<\Sym_n$.

\subsection{Intransitive actions}\label{sec:intransitive_actions}

\medskip

\begin{proposition} Let $2 \leqslant k < n-l < n$ and $n>3k+1$. If $l \geqslant 2k+1$ then $H= \langle r_n, r_{n-l}, r_k\rangle < \Sym_n$.

\end{proposition}

\begin{proof} Let $$\Phi = \{1, \ldots, k\} \cup r_n(\{1, \ldots, k\})=\{1, \ldots, k\} \cup \{n-k+1, \ldots, n\}.$$ Then $|\Phi|=2k$ and $\Phi$ is an $\langle r_n, r_k\rangle$-invariant set. Let $$\Delta=\{t \in \Omega : t \equiv s\pmod{l} \mbox{ for } s \in \Phi\}.$$ We claim that $\Delta$ is a proper $H$-invariant subset of $\Omega$ (i.\,e., $\Delta$  is a union of $H$-orbits). Indeed, let $x\equiv s\pmod{l}$, where $s\in \Phi$. We have $x = s + w\cdot l$ for some $w$ and $r_n(x)=n-x+1=n-s+1-w\cdot l \equiv n-s+1 \pmod{l}$, and $n-s+1=r_n(s) \in \Phi$. If $x \leqslant k$ then $r_k(x)=k+1-x=k-s+1-w\cdot l\equiv k-s+1 \pmod{l}$, and $ k-s+1=r_k(s)\in \Phi$; if $x>k$ then $r_k(x)=x$. If $x\leqslant n-l$ then $r_{n-l}(x)=n-x-l+1 \equiv r_n(x) \pmod{l}$, otherwise $r_{n-l}(x)=x$. Thus, $\Delta$ is a union of $H$-orbits. Moreover, since $|\Phi|<l$ we see that $\Delta$ is a proper subset of $\Omega$. Thus, $H< \Sym_n$.  
\end{proof}

We illustrate this Proposition~1 by the following example.

\begin{example}
Take $n=15$, $l=5$, and $k=2$. Then $\Omega=\{1,2,\ldots,15\}$, $2\leqslant k<n-l<n$, $n>3k+1$, and $l=2k+1$. Consider $H=\langle r_{15},r_{10},r_2\rangle$. Following the construction in the proof of Proposition~$1$, let
\[
\Phi
=
\{1,2\}\cup r_{15}(\{1,2\})
=
\{1,2,14,15\}.
\]
The set $\Phi$ is an $\langle r_{15},r_2\rangle$-orbit. Modulo $5$, the elements of \ $\Phi$ represent the residue classes $\{0,1,2,4\}$. Therefore
\begin{align*}
\Delta&=\{x\in\Omega:x\equiv0,1,2,\text{ or }4\pmod5\}=\\
&=\{1,2,4,5,6,7,9,10,11,12,14,15\}.
\end{align*}
Its complement is $\Omega\setminus\Delta=\{3,8,13\}$. Notice that $r_{15}(\Delta)=r_{10}(\Delta)=r_2(\Delta)=\Delta$, and so $\Delta$ is invariant under all three generators of $H$. Thus, $H$ is intransitive on $\Omega$ and consequently $\langle r_{15},r_{10},r_2\rangle<S_{15}$.
\end{example}

Not all examples of proper $H$-invariant subsets follow the construction from the proof of Proposition~1 as the following example illustrates.

\begin{example}
 Consider the prefix reversals of lengths $2$, $13$, and $18$ in $\Sym_{18}$ with $\Omega=\{1,2,\ldots,18\}$. In the notation of the proposition, we take $n=18, k=2, n-l=13$, and hence $l=5=2k+1$. Furthermore, let us consider the subset $\Delta=\{4,5,9,10,14,15\}\subsetneq\Omega$. Note that each of the three generators preserves $\Delta$ as a set. Indeed, we have:
 \[
 r_{18}(\Delta)
  =\{15,14,10,9,5,4\}
  =\Delta
 \quad
 \text{ and }
 \quad
 r_{13}(\Delta)
  =\{10,9,5,4,14,15\}
  =\Delta,
 \]
 where $r_{13}$ fixes $14$ and $15$. Moreover, $r_2$ fixes every
 element of $\Delta$, since all of its elements are greater than $2$. Consequently, $\Delta$ is invariant under $H=\langle r_{18},r_{13},r_2\rangle$.

 Equivalently, beginning with the identity permutation and applying any
 sequence of these three prefix reversals, the values $4,5,9,10,14,15$ must collectively occupy the positions
 $4,5,9,10,14,15$. Thus, every resulting permutation has the following form in one-line notation:
 \[
  [\star \,\star\,\star\,a\,b\,\star\,\star\,\star\,c\,d\,\star\,\star\,\star\,e\,f\,\star\,\star\,\star], 
 \quad \text{ with }\quad 
 \{a,b,c,d,e,f\}=\{4,5,9,10,14,15\},
 \]
 where the entries denoted by $\star$ are the elements of
 $\Omega\setminus\Delta$. Since $\Delta$ is a nonempty proper
 $H$-invariant subset of $\Omega$, the action of $H$ is intransitive. Therefore, $H<\Sym_{18}$.
\end{example}

Similarly, we have the following proposition.

\begin{proposition}\label{Divisors} Let $H=\langle r_n, r_{n-l},r_k\rangle$ for $n \geqslant 6$ and $l \geqslant k+1$. If $l$ divides $n+1$, then  $$\Delta =\{t : l \mbox{ divides } t\}$$ is a proper $H$-invariant subset of $\Omega$ {\rm(}i.\,e., $\Delta$ is a union of $H$-orbits{\rm)}, therefore $H < \Sym_n$.
\end{proposition}

\begin{proof} 
If $l$ divides $n+1$ and $l$ divides $i$ then $r_k(i)=i$, $l$ divides $r_n(i)=n-i+1$ and $l$ divides $r_{n-l}(i)=n-l-i+1$. So, $H$ preserves $\Delta$.
\end{proof}

We illustrate this proposition by the following example.

\begin{example}
Consider the prefix reversals of lengths $2$, $8$, and $11$ in $\Sym_{11}$, with $\Omega=\{1,2,\ldots,11\}$. In the notation of the proposition, we have $n=11$, $k=2$, $n-l=8$, so $l=3$. Thus, $l=k+1$ and $l\mid n+1$. The corresponding subset is as follows:
\[
\Delta=\{t\in\Omega:3\mid t\}=\{3,6,9\}.
\]
Each of the three generators preserves $\Delta$ as a set. Indeed, we have:
\[
r_{11}(\Delta)
 =\{9,6,3\}
 =\Delta
\quad
\text{ and }
\quad
r_8(\Delta)
 =\{6,3,9\}
 =\Delta.
\]
Finally, $r_2$ fixes each element of $\Delta$ since all three elements are greater than $2$. Hence, $\Delta$ is invariant under
$H=\langle r_{11},r_8,r_2\rangle.$

Equivalently, beginning with the identity permutation and applying any sequence of these prefix reversals, the values $3$, $6$, and $9$ must collectively occupy positions $3$, $6$, and $9$. Thus, every resulting permutation has the following form:
\[
[\star\,\star\,a\,\star\,\star\,b\,\star\,\star\,c\,\star\,\star],
\quad\text{ where }\quad
\{a,b,c\}=\{3,6,9\},
\]
where the entries denoted by $\star$ are the elements of
$\Omega\setminus\{3,6,9\}$. Since $\Delta$ is a nonempty proper $H$-invariant subset of $\Omega$, the action of $H$ is intransitive, and therefore $H<\Sym_{11}$.
\end{example}

Since $r_k$ for $k\in \{2,3\}$ is a transposition, by Lemma~\ref{WieLemma} we have $H = \Sym_n$ if and only if $H$ is primitive on $\Omega$. The following two lemmas are useful to establish which actions of the groups of the form $\langle r_n,r_m,r_k\rangle$, with $k<m<n$ and $k\in\{2,3\}$ acting naturally on $\Omega$, are intransitive. 

\begin{lemma}\label{lem:n-4_2}
    If the set $\{r_n,r_m,r_2\}$, where $n\geqslant 7$ and $2<m<n$, generates $\Sym_n$ then $m\geqslant n-4$. Furthermore, if $n$ is odd then $m\geqslant n-3$.
\end{lemma}

\begin{proof}
Let $d=n-m$ and consider the following union of residue classes modulo $d$:
\[
\Delta
=
\left\{
x\in\Omega:
x\equiv 1,2,n-1,\text{ or }n\pmod d
\right\}.
\]
We claim that $\Delta$ is invariant under $H=\langle r_n,r_m,r_2\rangle$. First, notice that $r_2$ swaps $1$ and $2$ and fixes every element greater
than $2$ so it preserves $\Delta$. Now let $t\in\{m,n\}$. If $x>t$ then $r_t(x)=x$ and there is nothing to prove. If $x\leqslant t$ then $r_t(x)=t+1-x$. Moreover, since $m\equiv n\pmod d$ then we have:
\[
r_t(x)\equiv n+1-x\pmod d.
\]
The involution $s\mapsto n+1-s\pmod d$ permutes the four residue classes defining $\Delta$. Thus, $r_m$ and $r_n$ preserve $\Delta$ proving that $\Delta$ is $H$-invariant.

If $d\geqslant 5$ then $\Delta$ consists of at most four of the $d$ residue classes modulo $d$. Since $d<n$ then every residue class modulo $d$ occurs among the elements of $\Omega$. Hence, $\Delta$ is a proper nonempty $H$-invariant subset of $\Omega$, so $H$ is intransitive. Therefore, if $H=\Sym_n$ then $d \leqslant 4$ or equivalently $m\geqslant n-4$.

It remains to consider the case $d=4$. In this case the residue classes defining $\Delta$ are $1, 2, n-1, n\pmod4$. 
If $n\equiv1\pmod4$ then they are $0,1,2\pmod4$, while if $n\equiv3\pmod4$ they are $1,2,3\pmod4.$ Thus, when $n$ is odd and $d=4$, the set $\Delta$ is again a proper $H$-invariant subset, so $H$ is intransitive. Consequently, if $n$ is odd and $H=\Sym_n$ then $d\leqslant 3$ which is equivalent to $m\geqslant n-3$.
\end{proof}

\begin{lemma}\label{rem:6k+1/5}
Suppose that $n\geqslant 6$,  $n\equiv 2\pmod 3$ and $m=n-3$. Then $\langle r_n,r_{n-3},r_2\rangle<\Sym_n$.
\end{lemma}
\begin{proof} Notice that $m\equiv n\equiv2\pmod3$. Let $\Delta=\{x\in\Omega:x\not\equiv0\pmod3\}$. We claim that $\Delta$ is invariant under $H=\langle r_n,r_{n-3},r_2\rangle$. Indeed, let $t\in\{n,n-3\}$. If $x\leqslant t$ then $r_t(x)=t+1-x\equiv -x\pmod3$. So, $r_t$ swaps the residue classes $1$ and $2$ modulo $3$. If $x>t$ then $r_t(x)=x$. Hence, both $r_n$ and $r_{n-3}$ preserve $\Delta$. The reversal $r_2$ also preserves $\Delta$. Since $3\notin\Delta$ then the set $\Delta$ is a proper subset of $\Omega$. Thus, $H$ is intransitive, and consequently $\langle r_n,r_{n-3},r_2\rangle<\Sym_n$.
\end{proof}

Similarly, we have the following lemma.

\begin{lemma}\label{lem:n-4_3}
    If the set $\{r_n,r_m,r_3\}$, where $3<m<n$ and $n\geqslant 8$, generates $\Sym_n$ then $m\geqslant n-4$. Furthermore, if $n$ is odd then $m\geqslant n-3$.
\end{lemma}

\begin{proof}
The proof is analogous to that of Lemma~\ref{lem:n-4_2}. As before, take $d=n-m$ and consider $\Delta = \{x\in[n]:x\equiv1,3,n-2,n\pmod d\}$. As above, we prove that the set $\Delta$ is invariant under $\langle r_n,r_m,r_3\rangle$.
\end{proof}

\section{Main results}\label{sec:results}

We now re-state and prove the main theorems. 

\subsection{Case $k=2$}\label{k23}

\medskip
\begin{theorem}\label{thm:2}
    The set $\{r_n,r_m,r_2\}$, where $2<m<n$ and $n\geqslant 4$, generates $\Sym_n$ if and only if
    \begin{enumerate}[(i)]
        \item $n$ is even and $m=n-1$, or 
        \item $n$ is odd and 
        \begin{enumerate}
            \item $m\in\{n-1,n-2\}$ if $n\equiv2\pmod 3$ or if $n\equiv0\pmod 3$, or
            \item $m\in\{n-3,n-2,n-1\}$ if $n\equiv 1\pmod 3$. 
        \end{enumerate}
    \end{enumerate}
\end{theorem}

\begin{proof} The cases with $n\leqslant 7$ are easily verified to hold. So, we assume that $n\geqslant 7$. The case $m=n-2$ was studied in~\cite{Gun22}, where it was proved that $H=\langle r_n, r_{n-2}, r_2\rangle=\Sym_n$ for any odd $n \geqslant 5$, and the case $m=n-1$ was considered in~\cite{BS03} where it was shown that $H=\langle r_n, r_{n-1}, r_2\rangle=\Sym_n$ for any  $n \geqslant 4$. Thus, we can assume that $m \leqslant n-3$.

If $n$ is even then $m\geqslant n-4$ by Lemma~\ref{lem:n-4_2}. Moreover, not all three prefix reversals can have an even index if they generate $\Sym_n$ (see Observation (5) in the Introduction); therefore, the only possible value remaining for $m$ is $n-3$. If $n$ is odd then by Lemma~\ref{lem:n-4_2} we have $m \geqslant n-3$. Thus, all that remains is the case $m=n-3$.

If $n\equiv0\pmod 3$ then the set 
$$\mathcal{B}_a=\{3a+b:1\leqslant b\leqslant 3\} \mbox{ with } 0\leqslant a\leqslant n/3-1$$ 
forms a non-trivial block system of $\langle r_n, r_{n-3}, r_2\rangle$ acting on $\Omega$. Hence, $\langle r_n, r_{n-3}, r_2\rangle<\Sym_n$ by Lemma~\ref{WieLemma}.

If $n\equiv 2\pmod 3$ then Lemma~\ref{rem:6k+1/5} gives that $\langle r_n, r_{n-3}, r_2\rangle<\Sym_n$.

If $n\equiv 4\pmod 6$ then the set $$\mathcal{B}_b=\{6a+b, 6a+b+1, 6a+b+2 :a\geqslant 0\} \mbox{ with }b \in \{0,3\}$$ forms a non-trivial block system of $\langle r_n, r_{n-3}, r_2\rangle$ acting on $\Omega$. Therefore, $\langle r_n, r_{n-3}, r_2\rangle< \Sym_n$. 

If $n\equiv 1\pmod 6$ then $r_nr_{n-3}r_2$ is the $n$-cycle as follows:
    \[
    \sigma=(1\ 5\ 8\ 11 \ \ldots \ (6s-1)\ 3 \ 6 \ 9\ 12 \ \ldots \ 6s \ 2 \ 4 \ 7 \ 10 \ 13 \ \ldots (6s+1)),  
    \]
    where $s=\lfloor n/6\rfloor$. To see this, notice that $r_nr_{n-3}r_2(1)=5$, $r_nr_{n-3}r_2(2)=4$. Since $r_2$ fixes every $i$ with $3\leqslant i\leqslant n-3$ then $r_nr_{n-3}r_2(i)=i+3$ for $3\leqslant i\leqslant n-3$. Moreover, $r_nr_{n-3}r_2(n-2)=r_n(n-2)=3$, $r_nr_{n-3}r_2(n-1)=r_n(n-1)=2$, and $r_nr_{n-3}r_2(n)=r_n(n)=1$. Thus, the following holds:
    \[
    r_nr_{n-3}r_2(i)=
    \begin{cases}
       5,&\text{ if $i=1$;} \\
       4,&\text{ if $i=2$;} \\
       i+3,&\text{ if $3\leqslant i\leqslant n-3$;}\\
       n+1-i,&\text{ if $n-2\leqslant i\leqslant n$.}
    \end{cases}
    \] From this, starting with $1$ we have:
    \[1\rightarrow5\rightarrow8\rightarrow 11\rightarrow\cdots\rightarrow 6s-1=n-2,\] and starting with $n-2$ one can continue as follows: \[n-2\rightarrow 3\rightarrow6\rightarrow\cdots\rightarrow6s=n-1,\] such that starting with $n-1$ we get the following chain: \[n-1\rightarrow 2\rightarrow 4\rightarrow7\rightarrow10\rightarrow\cdots\rightarrow 6s+1=n.\] Thus, finally we have $n\rightarrow1$ which gives the sought permutation $\sigma$. 

From the previous computation, notice that $(r_nr_{n-3}r_2)^{2s}(1)=3$, and so we have $(r_nr_{n-3}r_2)^{4s}(1)=(r_nr_{n-3}r_2)^{2s}(3)=2$. Moreover, $\gcd(4s,n)=\gcd(4s,6s+1)=1$. Since $r_2=(1\ 2)$ it follows that $\langle r_n,r_{n-3},r_2\rangle=\Sym_n$ by Lemma~\ref{lem:janusz}. \end{proof}

\subsection{Case $m=n-1$}\label{m=n-1}

\medskip

\begin{theorem}\label{thm:n-1}
	Let $H=\langle r_n, r_{n-1}, r_k\rangle$, where $2\leqslant k\leqslant n-2$ and $n\geqslant 4$.
	Then $H=\Sym_n$ if and only if
	\begin{enumerate}
		\item[$(i)$] $n$ is even and $k$ is even, or 
		\item[$(ii)$] $n \equiv 3 \pmod{4}$ and  $2\leqslant k\leqslant n-2$, or
		\item[$(iii)$] $n \equiv 1 \pmod{4}$ and $k \equiv 2, 3 \pmod{4}$. 
	\end{enumerate}
\end{theorem}

\begin{proof} 

Consider the following two permutations:
\begin{align*}
	R &= (1\ 2 \dots n-1\ n)=r_{n}r_{n-1} = L^{-1},\\
	L &= (n\ n-1\dots 2\ 1) =r_{n-1}r_{n} = R^{-1}.
\end{align*}

By direct computations, we obtain the following 3-cycle:
\begin{align*}
	\sigma &= r_k R r_k R r_k L r_k L = (1 \ k\ k+2).
\end{align*}

We have:
$$
d = \gcd(k-1, k+1, n) = \gcd(2, k+1, n) = 
\begin{cases}
	1, & \text{if $n$ is odd}; \\[6pt]
	1, & \text{if $n$ is even and $k$ is even}; \\[6pt]
	2, & \text{if $n$ is even and $k$ is odd}.
\end{cases}
$$

Let $H_1 = \langle L, \sigma \rangle \leqslant H$.
By Lemma~\ref{lem:Gen_(3,n)}, $H_1 \geqslant \mathrm{Alt}_n$ if and only if $d = 1$. 

Let us show that if $d=2$ then $H_1$ and $H$ are imprimitive.

If $n$ is odd then $H_1 \geqslant \mathrm{Alt}_n$. If $n \equiv 3 \pmod{4}$ then $r_n \in \Sym_n \setminus \mathrm{Alt}_n$ and, therefore, $\mathrm{Alt}_n < H = \Sym_n$. If $n \equiv 1 \pmod{4}$ then $r_n, r_{n-1} \in \mathrm{Alt}_n$, and $r_k \not \in \mathrm{Alt}_n$ if and only if $k \equiv 2, 3 \pmod{4}$.

If both $n$ and $k$ are even then $H \geqslant H_1=\Sym_n$ by Lemma~\ref{lem:Gen_(3,n)}.

It remains to show that in the case of even $n$ and odd $k$ the group $H$ is imprimitive. Let $D = \langle L, r_n \rangle \leqslant H$. Then $D$ is isomorphic to the dihedral group $D_{2n}$ of order $2n$, and we claim that the following two sets: 
$$\mathcal{B}_1 = \{t \mid t \equiv 1\pmod{2}\}  \mbox{ and } \mathcal{B}_2 = \{t \mid t \equiv 0\pmod{2}\}$$ 
are blocks of imprimitivity of $D$. Indeed, $L(\mathcal{B}_i)=\mathcal{B}_j$ and $r_n(\mathcal{B}_i)=\mathcal{B}_j$ for $j\in \{1,2\}\setminus\{i\}$. Now we claim that $r_k$ and $r_{n-1}$ also preserve blocks $\mathcal{B}_1$ and $\mathcal{B}_2$. Indeed, since both $k$ and $n-1$ are odd then $r_k(\mathcal{B}_i)=\mathcal{B}_i$ and $r_{n-1}(\mathcal{B}_i)=\mathcal{B}_i$ for each $i \in \{1,2\}$.  Thus, $\{\mathcal{B}_1, \mathcal{B}_2\}$ is a non-trivial block system of $H$ acting on $\Omega$, and so $H<\Sym_n$.
\end{proof}

\subsection{Case $k=3$}\label{k3}

\medskip

\begin{theorem}\label{thm:3}
The set $\{r_n,r_m,r_3\}$, with $2\leqslant m<n$, $m \not =3$ and $n\geqslant 5$, generates $\Sym_n$ if and only if one of the following conditions holds{\rm:}
    \begin{enumerate}[$(i)$]
        \item $n$ is even and $m=n-2${\rm;}
        \item $m\in\{n-3,n-1\}$ and $n\equiv 1,5\pmod 6${\rm;}
        \item $m=n-1$ if $n\equiv3\pmod 6$.
    \end{enumerate}
\end{theorem}

\begin{proof} The case $\{r_n, r_3, r_2\}$ is covered by Theorem~1, and the only triple which generates $\Sym_n$ in this case is $\{r_5, r_3, r_2\}$ when $n=5$. Thus, one can assume that $r_2$ is not included in a triple of prefix reversals. The cases $n\leqslant 8$ are easily verified by direct computation, therefore we can assume that $n \geqslant 9$. 

By Lemma~\ref{lem:n-4_3}, it follows that $m\geqslant n-4$. 
Moreover, it was shown in~\cite{BS03} that $\langle r_n, r_{n-2}, r_3\rangle=\Sym_n$ for any even $n \geqslant 6$ and $\langle r_n,r_{n-1},r_3\rangle =\Sym_n$ for any odd $n\geqslant5$. These establish all the sufficient conditions except $m=n-3$. 

Let us consider the case $m=n-4$ and $n$ even, for which we show that $\langle r_n,r_{n-4},r_3\rangle<\Sym_n$. Notice that $n$ must be congruent to $0$ or $2$ modulo $4$. If $n\equiv 0\pmod 4$ then the set \[\mathcal{B}_a=\{4a+1,4a+2,4a+3,4a+4\}\text{ with }0\leqslant a\leqslant n/4-1\] forms a non-trivial block system of $\langle r_n,r_{n-4},r_3\rangle$ acting on $\Omega$. 

Additionally, if $n\equiv6\pmod 8$ then the set 
\[\mathcal{B}_1=\{8a+b:a\geqslant0\text { and }0\leqslant b\leqslant 3\}\] along with the set
\[\mathcal{B}_2=\{8a+b:a\geqslant0\text{ and } 4\leqslant b\leqslant 7\}\] form a non-trivial block system of $\langle r_n,r_{n-4},r_3\rangle$ acting on $\Omega$. 

Finally, if $n\equiv 2\pmod8$ then the set
\[\mathcal{B}_1=\{8a+b:a\geqslant0 \text{ and }b\in\{1,3,4,6\}\}\] along with the set
\[\mathcal{B}_2=\{8a+b:a\geqslant0\text{ and }b\in\{0,2,5,7\}\}\] form a non-trivial block system of $\langle r_n,r_{n-4},r_3\rangle$ acting on $\Omega$. Thus, $\langle r_n,r_{n-4},r_3\rangle<\Sym_n$ by Lemma~\ref{WieLemma}. 

\smallskip

Furthermore, if $n$ is even then both $n-1$ and $n-3$ are odd. Therefore, $\langle r_n,r_{n-1},r_3\rangle <\Sym_n$ by Theorem~\ref{thm:n-1}, where in its proof it is shown that $\langle r_n,r_{n-1},r_3\rangle$ is imprimitive. For the case $\{r_n,r_{n-3},r_3\} $, consider the  partition $\mathcal{B}_{e},\mathcal{B}_{o}$ of $\Omega$, where $\mathcal{B}_e=\{2,4,\ldots,n\}$ and $\mathcal{B}_o=\{1,3,\ldots,n-1\}$. It is easy to see that $\{\mathcal{B}_o,\mathcal{B}_e\}$ is a nontrivial block system. Indeed, since $r_3$ and $r_{n-3}$ are odd-indexed they preserve the blocks; furthermore, $r_n(\mathcal{B}_e)=\mathcal{B}_o$. Thus, $\langle r_n,r_{n-3},r_3\rangle<\Sym_n$. 

\smallskip

Now let us suppose that $n$ is odd. By Lemma~\ref{lem:n-4_3} it follows that $m\geqslant n-3$. If $m=n-2$ then $\{r_n, r_{n-2}, r_3\}$ cannot generate $\Sym_n$ since all the numbers $3$, $n-2$ and $n$ are odd (see Observation (3) in the Introduction). Thus, the only remaining values to consider are $m=n-1$ and $m=n-3$. Moreover, if $n\equiv 3\pmod 6$ then $\{r_n,r_{n-3},r_3\}$ cannot generate $\Sym_n$ since $\gcd(n,n-3,3)=3$ (see Observation (5) in the Introduction). 

Again, we know that $\langle r_n,r_{n-1},r_3\rangle=\Sym_n$ for all odd $n\geqslant 5$~\cite{BS03}. We are now only left with the cases $n\equiv 1,5\pmod 6$ and $m=n-3$, and we now proceed to show that all of these generate $\Sym_n$.

If $n\equiv 1\pmod 6$ then $r_nr_{n-3}r_{3}$ is the $n$-cycle as follows (the proof is similar to that in Theorem~1):
    \[
    \sigma_1=(1\ 6\ 9\ 12 \ \ldots \ 6s\ 2\ 5\ 8 \ 11 \ \ldots \ (6s-1)\ 3\ 4\ 7\ 10\ \ldots\ (6s+1)),
    \] where $s=\lfloor n/6\rfloor$. To see this, we observe that $r_nr_{n-3}r_3(1)=6, r_nr_{n-3}r_3(2)=5$, and $r_nr_{n-3}r_3(3)=4$ using right-to-left composition. Since $r_3$ fixes every $i$ with $4\leqslant i\leqslant n-3$ then we have $r_nr_{n-3}r_3(i) = r_nr_{n-3}(i)=r_n(n-2-i)=i+3$ for $4\leqslant i\leqslant n-3$. Moreover, $r_3$ and $r_{n-3}$ both fix $\{n-2,n-1,n\}$, and hence $r_nr_{n-3}r_3(n-2)=r_n(n-2)=3$, $r_nr_{n-3}r_3(n-1)=r_n(n-1)=2$, and $r_nr_{n-3}r_3(n)=r_n(n)=1$. Therefore, we have:
\[
r_nr_{n-3}r_3(i)=
\begin{cases}
6, & \text{if }i=1;\\
5, & \text{if }i=2;\\
4, & \text{if }i=3;\\
i+3, & \text{if }4\leqslant i\leqslant n-3;\\
n+1-i, & \text{if }n-2\leqslant i\leqslant n.
\end{cases}
\]

Starting with $1$, we obtain the following chain: \[1\rightarrow6\rightarrow9\rightarrow12
\rightarrow\cdots\rightarrow6s=n-1,\] where the cycle can be continued as follows: \[n-1\rightarrow2\rightarrow5\rightarrow8
\rightarrow11\rightarrow\cdots
\rightarrow6s-1=n-2,\] and further as follows:
\[n-2\rightarrow3\rightarrow4\rightarrow7
\rightarrow10\rightarrow\cdots
\rightarrow6s+1=n.\]
Thus, finally have $n\rightarrow 1$ which gives us the sought $n$-cycle $\sigma_1$.

Notice that $\sigma_1^{4s}(1)=3$ and $\gcd(4s,n)=\gcd(4s,6s+1)=1$. Thus,  since $r_3=(1\ 3)$ then by Lemma~\ref{lem:janusz} the set $\{r_n,r_{n-3},r_3\}$ generates $\Sym_n$.

Similarly, if $n\equiv 5\pmod 6$ then $r_nr_{n-3}r_{3}$ is the following $n$-cycle:
    \[
    \sigma_2=(1\ 6\ 9\ 12 \ \ldots \ (6s+3)\ 3\ 4 \ 7 \ 10 \ \ldots \ (6s+4)\ 2\ 5\ 8\ 11 \ \ldots \ (6s+5)),
    \] where $s=\lfloor{n/6}\rfloor$. To see this, we observe that $r_nr_{n-3}r_3(1)=6, r_nr_{n-3}r_3(2)=5$, and $r_nr_{n-3}r_3(3)=4$ using right-to-left composition. Since $r_3$ fixes every $i$ with $4\leqslant i\leqslant n-3$ then we have $r_nr_{n-3}r_3(i)=r_nr_{n-3}(i)=r_n(n-2-i)=i+3$ for $4\leqslant i\leqslant n-3$. Moreover, $r_3$ and $r_{n-3}$ both fix $\{n-2,n-1,n\}$, and hence $r_nr_{n-3}r_3(n-2)=r_n(n-2)=3$, $r_nr_{n-3}r_3(n-1)=r_n(n-1)=2$, and $r_nr_{n-3}r_3(n)=r_n(n)=1$. Therefore, we have:
\[
r_nr_{n-3}r_3(i)=
\begin{cases}
6, & \text{if }i=1;\\
5, & \text{if }i=2;\\
4, & \text{if }i=3;\\
i+3, & \text{if }4\leqslant i\leqslant n-3,\\
n+1-i, & \text{if }n-2\leqslant i\leqslant n.
\end{cases}
\]

Starting with \(1\), we obtain the following chain:
\[
1\rightarrow6\rightarrow9\rightarrow12
\rightarrow\cdots\rightarrow6s+3=n-2.
\]
Since \(n-2\rightarrow3\) the cycle can be continued as follows:
\[
n-2\rightarrow3\rightarrow4\rightarrow7
\rightarrow10\rightarrow\cdots
\rightarrow6s+4=n-1,
\]
and since \(n-1\rightarrow2\) the cycle can be continued as follows: 
\[
n-1\rightarrow2\rightarrow5\rightarrow8
\rightarrow11\rightarrow\cdots
\rightarrow6s+5=n.
\]
Finally, $n\rightarrow1$ which gives the sought $n$-cycle $\sigma_2$.

As before, $\sigma_2^{2s+1}(1)=3$ and $\gcd(2s+1,n)=\gcd(2s+1,6s+5)=1$. Hence, $\{r_n,r_{n-3},r_3\}$ generates $\Sym_n$ by Lemma~\ref{lem:janusz}. This completes the proof. 
\end{proof}

\subsection{Case $m=n-2$}\label{m=n-2}

\medskip

\begin{theorem}\label{thm:n-2}
Let $H=\langle r_n, r_{n-2}, r_k\rangle$, with $2\leqslant k < n-2$ and $n\geqslant 5$. Then $H=\Sym_n$ if and only if $n$ and $k$ are of different parities. 
\end{theorem}

\begin{proof} The cases $k \in \{2, 3\}$ were already handled in Theorems~1 and~3, respectively. Therefore, one can assume that $k>3$. We first establish the sufficient statement. Suppose $n$ is even and $k$ is odd. Then $H$ contains the following $n$-cycle:
$$\sigma = r_{n-2}r_{n}=(n\ n-2\ n-4\dots 2\ n-1\ n-3 \dots \ 3 \  1).$$ 
To see this, observe that the following holds:
\[
\sigma=r_{n-2}r_n(i)=
\begin{cases}
	n, & \text{if }i=1;\\
	n-1, & \text{if }i=2;\\
	i-2, & \text{if }3\leqslant i\leqslant n.\\
\end{cases}
\]
Starting with \(n\) we obtain the chain:
\[
n\rightarrow n-2\rightarrow n-4 \rightarrow\cdots\rightarrow 2.
\]
Since \(2\rightarrow n-1\) the cycle can be continued as follows:
\[
n-1\rightarrow n-3\rightarrow n-5\rightarrow 3 \rightarrow\cdots \rightarrow 1.
\]
Finally, we have \(1\rightarrow n\) which gives the sought $n$-cycle $\sigma$.

Since $k$ is odd and $k > 3$ then for $k=2m+3$ let us consider the following permutation:
\[
\rho  = (r_{n-2}r_n)^{\frac{k-3}{2}}r_k = \sigma^{m}r_k.
\] 
For $i \geqslant 2m+1=k-2$, applying $\sigma^{m}$ simply subtracts $2$, $m$ times; thus, we have $\rho(1) = \sigma^{m}(2m + 3) = 3$, $\rho(2) = \sigma^{m}(2m + 2) = 2$ and $\rho(3) = \sigma^{m}(2m + 1) = 1$.
Therefore, the restriction of $\rho$ to the set $\{1,2,3\}$ satisfies the following: $$\rho\vert_{\{1,2,3\}} = (1\ 3).$$

Let $4 \leqslant i \leqslant n$. If $i \leqslant k$ then we have: 
\[
\rho(i) = \sigma^{m}r_k(i) = \sigma^{m}(k+1-i) = \sigma^{m}(2m+4-i),
\]
and
\[
\sigma^{m}(2x) = n-1-2(m-x), \; \text{for } 1 \leqslant x \leqslant m;
\]
\[
\sigma^{m}(2x+1) = n-2(m-x-1), \; \text{for } 0 \leqslant x \leqslant m-1.
\]

If $i$ is even then $2m+4-i$ is even and the following holds:
\[
\sigma^{m}(2m+4-i) = n-1 - 2m + 2m+4-i = n + 3 - i.
\]

If $i$ is odd then $2m+4-i$ is odd and the following holds:
\[
\sigma^{m}(2m+4-i) = n + 5 - i.     
\]

If $i > k = 2m+3$ then $\rho(i) = \sigma^{m}(i) = i - 2m$. Thus, we have:
\[
\rho(i)=
\begin{cases}
	3, & \text{if }i=1;\\
	2, & \text{if }i=2;\\
	1, & \text{if }i=3;\\
	n+3 - i, & \text{if }4\leqslant i\leqslant k \mbox{ and } i \text{ even};\\
	n+5 - i, & \text{if }4\leqslant i\leqslant k \mbox{ and } i \text{ odd};\\
	i - 2m, & \text{if }k < i\leqslant n.\\
\end{cases}
\]

Starting with \(4\) we obtain the chain:
\[
4 \rightarrow n-1 \rightarrow n-1-2m \rightarrow n-1-4m \rightarrow\cdots\rightarrow a_1, 
\]
where $a_1 \in \{4, \dots, k\}$ and $a_1 \equiv n-1 \pmod{2m}$. Since $n-1$ is odd then $a_1$ is odd and we have: 
\[
a_1 \rightarrow n + 5-a_1 \rightarrow n + 5-a_1 - 2m \rightarrow\cdots\rightarrow a_2,
\]
where $a_2 \in \{4, \dots, k\}$ and  $a_2 \equiv n + 5-a_1 \equiv  n + 5 -(n-1) \equiv 6 \pmod{2m}$. Thus, we have the following chain:
\[
4 \rightarrow \cdots \rightarrow 6 \rightarrow \cdots \rightarrow 8 \rightarrow \cdots \rightarrow k-1.
\]
Since $k-1 + 2 \equiv 2m+4 \equiv 4 \pmod{2m}$ we continue as follows: 
\[
k-1 \rightarrow \cdots \rightarrow 4.
\]
It remains to observe that each interval $a \rightarrow \cdots \rightarrow a+2$ contains exactly $2m$ terms and all of them belong to the same residue class modulo $2m$, namely, the class containing $n+3-a$. Thus, all odd values $\{5,7,\dots,k\}$ occur in the sequence, and consequently $\rho|_{\{4,5,\dots,n\}}$
is an $(n-3)$-cycle.

Hence, the following two equations hold: 
\begin{align*} (r_{n-2}r_n)^{\frac{k-3}{2}}r_k &= (1\ 3) \underbrace{(\dots \ldots \ldots \ldots \ldots)}_{\text{$(n-3)$-cycle}}, \\
((r_{n-2}r_n)^{\frac{k-3}{2}}r_k)^{n-3}&= (1\ 3). 
\end{align*}
Therefore, having the $n$-cycle $r_{n-2}r_n$ and the transposition $(1\ 3)$ that swaps two adjacent entries in this cycle, by Lemma~\ref{lem:janusz}, we have $H=\Sym_n$.

Now, let us assume that $n$ is odd and $k$ is even. Then we have: 
$$\sigma=r_{n-2} r_n=(n\ n-2\ n-4\dots 1)(n-1\ n-3 \dots 2).$$ 
Indeed, since $n$ is odd we have the following two chains:
\[
n \rightarrow n-2 \rightarrow n-4 \rightarrow\cdots\rightarrow 1 \rightarrow n
\]
and 
\[
n-1 \rightarrow n-3 \rightarrow n-5 \rightarrow\cdots\rightarrow 2 \rightarrow n-1.
\]

Since $k$ is even and $k > 2$, we can set $k=2m+2$ and consider the following permutation:
\[
\rho  = (r_{n-2}r_n)^{\frac{k-2}{2}}r_k = \sigma^{m}r_k
\]
with $\rho(1) = \sigma^{m}(2m+2) = 2$ and $\rho(2) = \sigma^{m}(2m+1) = 1$ which gives us: $$\rho\vert_{\{1,2\}} = (1\ 2).$$

Let $3 \leqslant i \leqslant n$. If $i \leqslant k$ then the following hold: 
\[
\rho(i) = \sigma^{m}r_k(i) = \sigma^{m}(k+1-i) = \sigma^{m}(2m+3-i) 
\]
and we have: 
\[
\sigma^{m}(2x) = n-1-2(m-x), \; \text{for } 1 \leqslant x \leqslant m
\]
and
\[
\sigma^{m}(2x+1) = n-2(m-x-1), \; \text{for } 0 \leqslant x \leqslant m-1.
\]

If $i$ is odd then $2m+3-i$ is even and we have: 
\[
\sigma^{m}(2m+3-i) = n-1 - 2m + (2m+3-i) = n + 2 - i.
\]
If $i$ is even then $2m+3-i$ is odd and we have: 
\[
\sigma^{m}(2m+4-i) = n + 4 - i.     
\]

If $i > k = 2m+2$ then $\rho(i) = \sigma^{m}(i) = i - 2m$. Hence, we have:
\[
\rho(i)=
\begin{cases}
	2, & \text{if }i=1;\\
	1, & \text{if }i=2;\\
	n+4 - i, & \text{if }3\leqslant i\leqslant k \mbox{ and } i \text{ even};\\
	n+2 - i, & \text{if }3\leqslant i\leqslant k \mbox{ and } i \text{ odd};\\
	i - 2m, & \text{if }k < i\leqslant n.\\
\end{cases}
\]

It remains to show that the restriction of $\rho$ to $\{3,4,\dots,n\}$ is a single $(n-2)$-cycle. Indeed, starting with $3$ we have the chain as follows:
\[
3 \rightarrow n-1 \rightarrow n-1-2m \rightarrow n-1-4m \rightarrow \cdots \rightarrow b_1, 
\]
where $b_1 \in \{3, \dots, k\}$ and $b_1 \equiv n-1 \pmod{2m}$. Since $n-1$ is even, we have $b_1$ is even and  
\[
b_1 \rightarrow n+4-b_1 \rightarrow n+4-b_1-2m \rightarrow n+4-b_1-4m \rightarrow \cdots \rightarrow b_2, 
\] 
where $b_2 \in \{3, \dots, k\}$ and $b_2 \equiv n + 4-b_1 \equiv  n + 4 -(n-1) \equiv 5 \pmod{2m}$. Thus, finally we have the following chain:
\[
3 \rightarrow \cdots \rightarrow 5 \rightarrow \cdots \rightarrow 7 \rightarrow \cdots \rightarrow k-1 \rightarrow 3.
\]

It remains to see that each interval $a \rightarrow \cdots \rightarrow a+2$ contains exactly $2m$ terms and all of them belong to the same residue class modulo $2m$, namely, the class containing $n+2-a$. Thus, all even values $\{4,6,\dots,k\}$ occur in the sequence, and consequently $\rho|_{\{3,4,\dots,n\}}$ is an $(n-2)$-cycle.

Thus, finally the following two equations hold: 
\begin{align*}
(r_{n-2}r_{n})^{\frac{k-2}{2}}r_k &= (1\ 2) \underbrace{(\dots \ldots \ldots \ldots \ldots)}_{(n-2)\text{-cycle}}, \\
((r_{n-2}r_{n})^{\frac{k-2}{2}}r_k)^{n-2} &=(1\ 2).
\end{align*}

Since the lengths of the two cycles in the permutation $r_{n-2}r_n$ are coprime as they differ exactly by $1$, it follows that for any odd $a$ and even $b$ there is $s \in \mathbb{Z}^+$ such that $(r_{n}r_{n-2})^{s} (1\ 2)(r_{n-2}r_{n})^{s} = (a\ b)$. Thus, $H$ contains the Coxeter generators of $\Sym_n$, and therefore $H = \Sym_n$.

For the necessary statement, notice that if both $n,k$ have the same parity, then all generators of $H$ have indices that are all of the same parity. Therefore, by Observations $(3)$ and $(5)$ in the Introduction, $H<\Sym_n$. This completes the proof. 
\end{proof}

\section{Conjectures and computational results}\label{sec:comp}

\subsection{Generating sets}\label{sec:generating_sets}

For a given $n$ and a set $\{r_n,r_m,r_k\}$, $2 \leqslant k<m<n$, it is natural to ask the following question: which fraction of pairs $\{r_m, r_k\}$ along with $r_n$ generate the entire $\mathrm{Sym}_n$? From computational results of a numerical experiment over all pairs $(m, k)$ for each $6 \leqslant n \leqslant 100$, we have the following conjecture. 

\begin{conjecture}
\label{conj:approx}
    Let $f(n)$ be the number of pairs $(m,k)$ such that $2 \leqslant k < m < n$ and $\langle r_n, r_m, r_k \rangle = \mathrm{Sym}_n$.
    Then 
    $$
     f(n) = \begin{cases}
        \frac{1}{13}n^2 + O(n) & n = 0 \pmod{4}; \\[6pt]
        \frac{1}{9}n^2 + O(n) & n = 1 \pmod{4}; \\[6pt]
        \frac{1}{10}n^2 + O(n) & n = 2 \pmod{4}; \\[6pt]
        \frac{1}{7}n^2 + O(n) & n = 3 \pmod{4}. \\[6pt]
    \end{cases}
    $$
\end{conjecture}

\begin{figure}[htbp]
\centering
\includegraphics[width=0.6\linewidth]{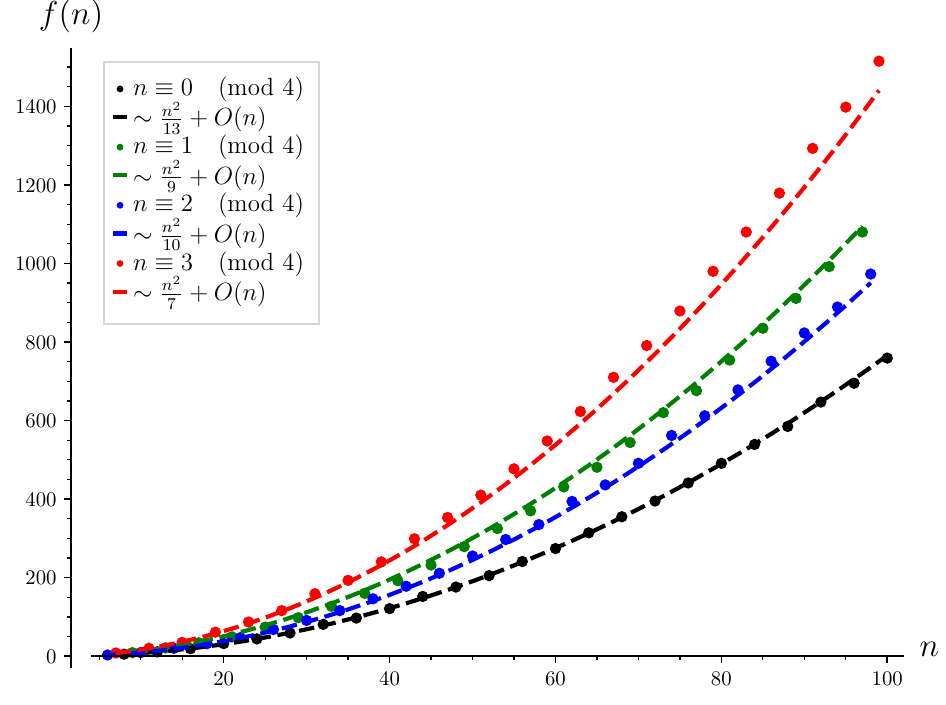}
 \caption{Visualization of Conjecture~\ref{conj:approx}. The figure depicts the number $f(n)$ of pairs $(m,k)$ so that $\langle r_n,r_m,r_k\rangle=\Sym_n$ for several values of $n$. Four colors are used representing the residue $r$ of $n$ modulo 4. The curves represent the best rational model $(x^2 + u x + v)/q$:\\
$r=0$ (black): $q=13$, $u=-0.5$, $v=5$,  $\mathrm{RMSE}=3.343$;\\
$r=1$ (green): $q=9$,  $u=5$,  $v=-45$,  $\mathrm{RMSE}=12.132$;\\
$r=2$ (blue): $q=10$, $u=-1$, $v=0$, $\mathrm{RMSE}=11.408$;\\
$r=3$ (red): $q=7$,  $u=3$,  $v=-14$, $\mathrm{RMSE}=38.677$.
}
\label{fig:conjecture-approx}
\end{figure}

Figure~\ref{fig:conjecture-approx} illustrates this conjecture. The plot displays the values of the function $f(n)$ for all $n \leqslant 100$.  Four distinct subsets are highlighted, each corresponding to a residue class of $n$ modulo $4$. For each subset, the best quadratic approximation with rational coefficients was computed using the Least Squares Method. For every fitted curve, the parameters of the quadratic function and the Root Mean Square Error (RMSE) metric are provided. 

In addition, there seems to be a connection between the sets of triples generating the whole group of degree $n$ modulo $4$ and modulo $2$. Based on results of computational experiments for $n \leqslant 100$ the following conjectures can be stated assuming that  $2\leqslant k < m < n$.

\begin{conjecture}
\label{conj:mod4}
If $\langle r_n, r_m,r_k\rangle=\mathrm{Sym}_{n}$ then $\langle r_{n+4},r_{m+4},r_{k+4}\rangle=\mathrm{Sym}_{n+4}$.
\end{conjecture}

\begin{conjecture}
\label{conj:mod2}
For any $n \equiv 0, 1 \pmod{4}$, if $\langle r_n,r_m,r_k\rangle=\mathrm{Sym}_{n}$ then $\langle r_{n+2},r_{m+2},r_{k+2}\rangle=\mathrm{Sym}_{n+2}$. 
\end{conjecture}

\begin{conjecture}
\label{conj:km_sum} If $m+k<n$ then $\langle r_n,r_m,r_k\rangle < \Sym_n$. Moreover, if $n$ is odd and $m+k<n-1$ then $\langle r_n,r_m,r_k\rangle < \Sym_n$.
\end{conjecture}

\begin{figure}[htbp]
    \centering
    \begin{subfigure}[b]{0.3\textwidth}
        \includegraphics[width=\textwidth]{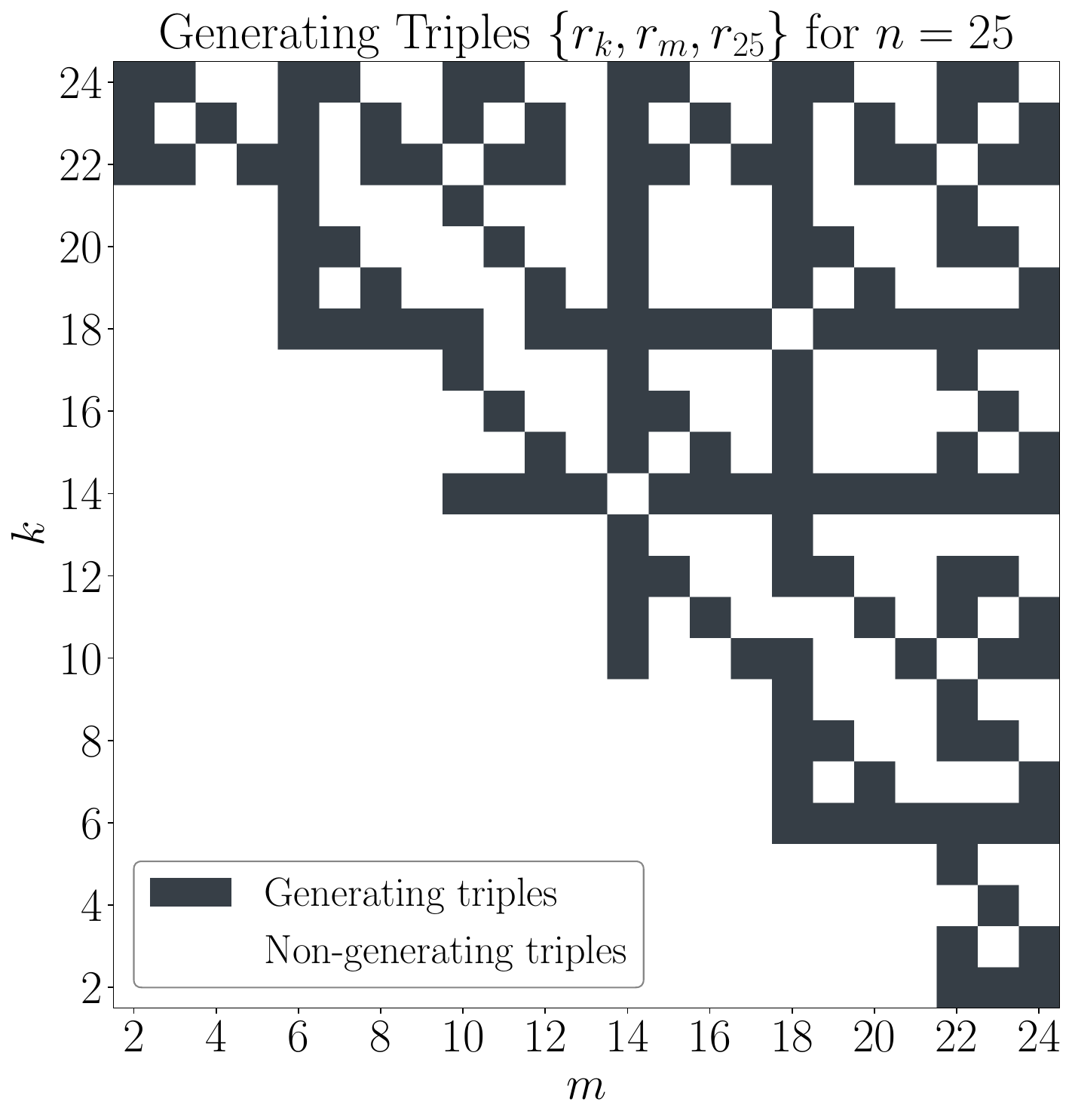}
        \label{fig:25}
    \end{subfigure}
    \hfill 
    \begin{subfigure}[b]{0.3\textwidth}
        \includegraphics[width=\textwidth]{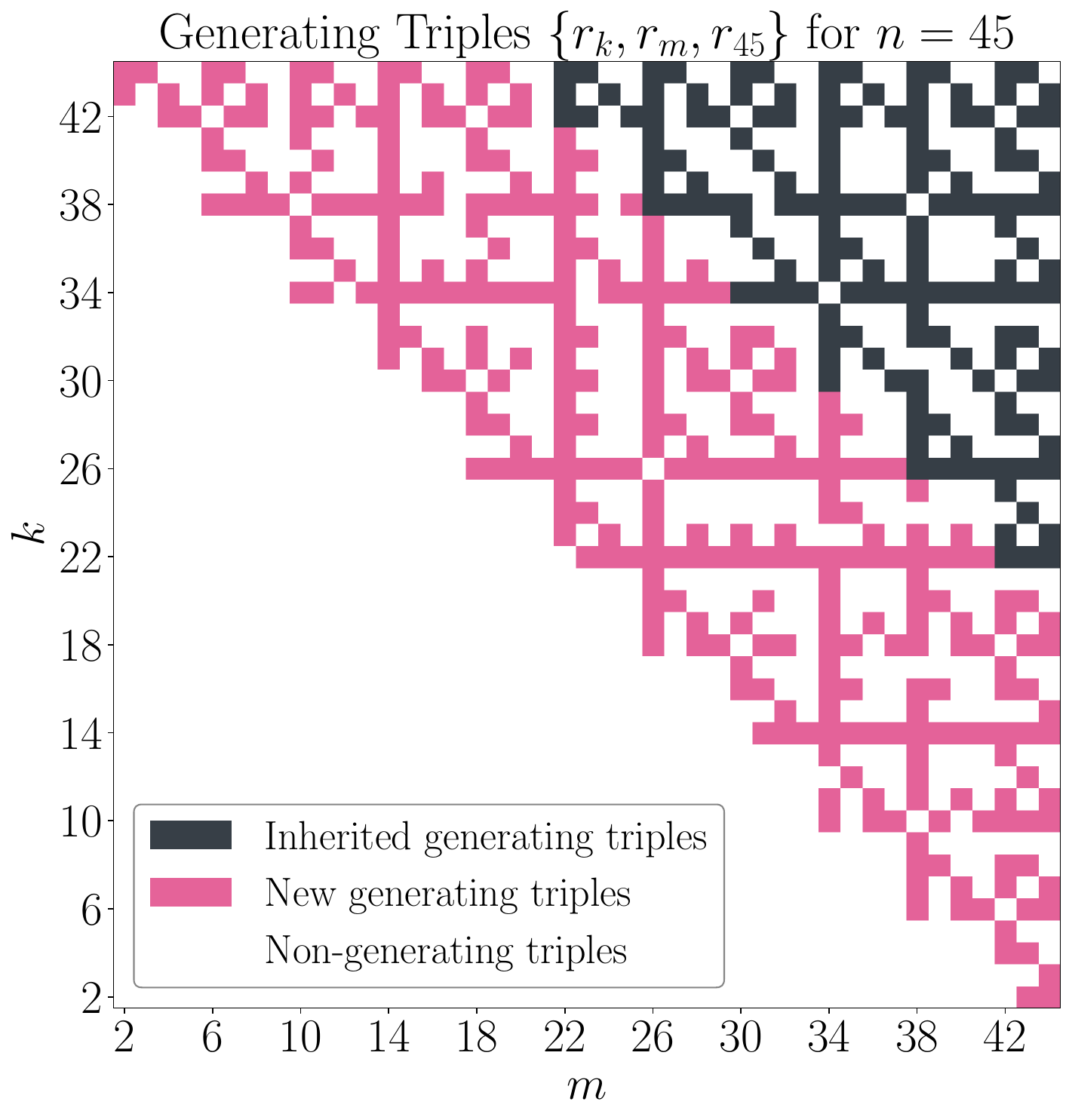}
        \label{fig:45}
    \end{subfigure}
    \hfill
    \begin{subfigure}[b]{0.3\textwidth}
        \includegraphics[width=\textwidth]{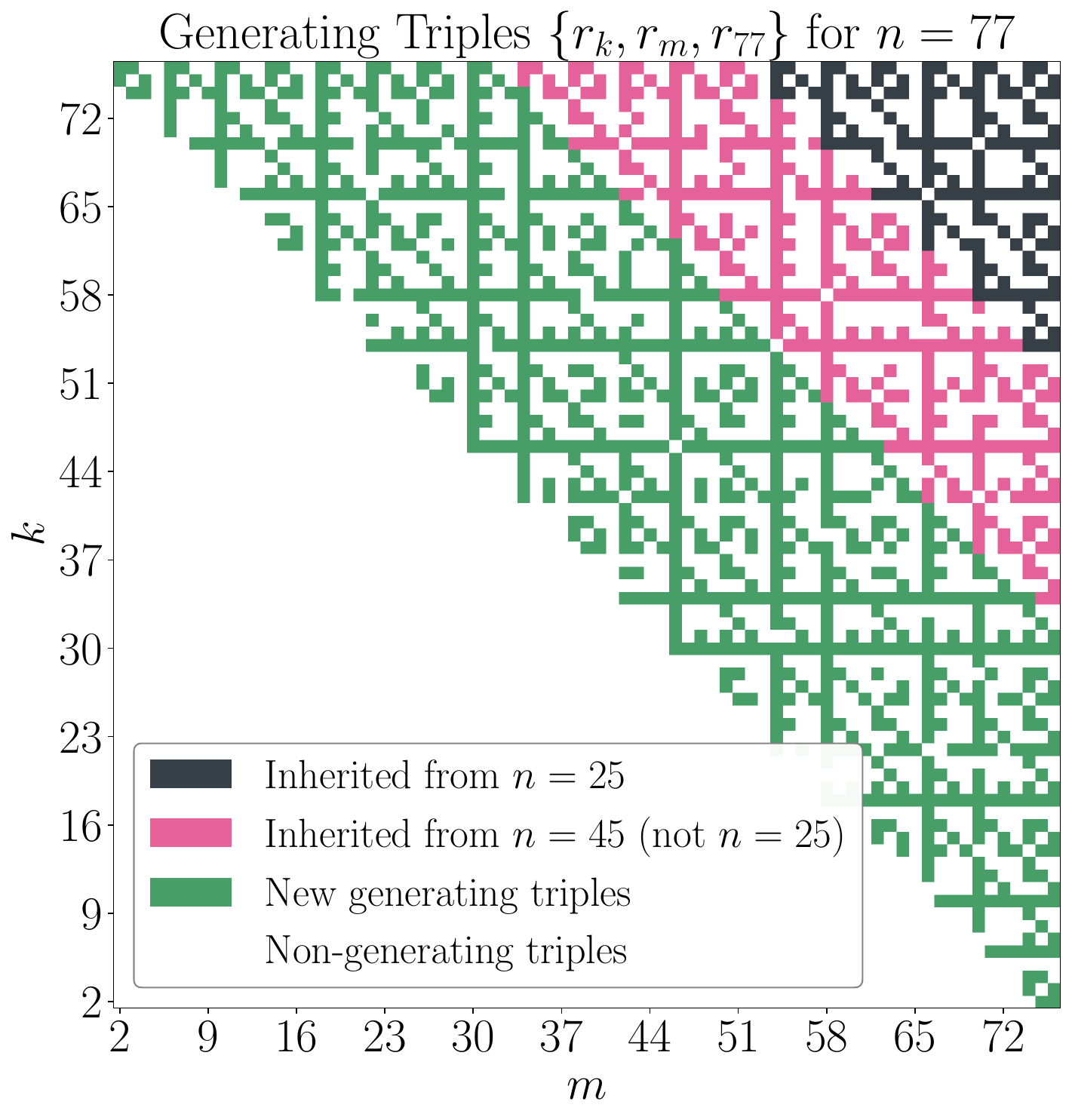}
        \label{fig:77}
    \end{subfigure}
    \caption{Visualization of Conjecture~\ref{conj:mod4} to illustrate that as $n$ grows, the previous triples are ``preserved'' $\pmod{4}$. The left-most figure depicts all pairs $(m,k)$ so that $\langle r_n,r_m,r_k\rangle=\Sym_n$. These pairs give rise to generating triples for $n=25+4\times 5$ by adding $4\times5$ to each index, and these are depicted in the middle figure (shown in black). There are other generating triples that are not inherited from the $n=45$ case (shown in pink). The right-most figure is built similarly by adding $4\times 8$ to each index for the triples for $n=45$. We see again that the previous pairs in black and pink are preserved for $n=77$ by adding $4\times8$ to each index.} 
    \label{fig:mod4}
\end{figure}

\begin{figure}[htbp]
    \centering
    \begin{subfigure}[b]{0.45\textwidth}
        \includegraphics[width=\textwidth]{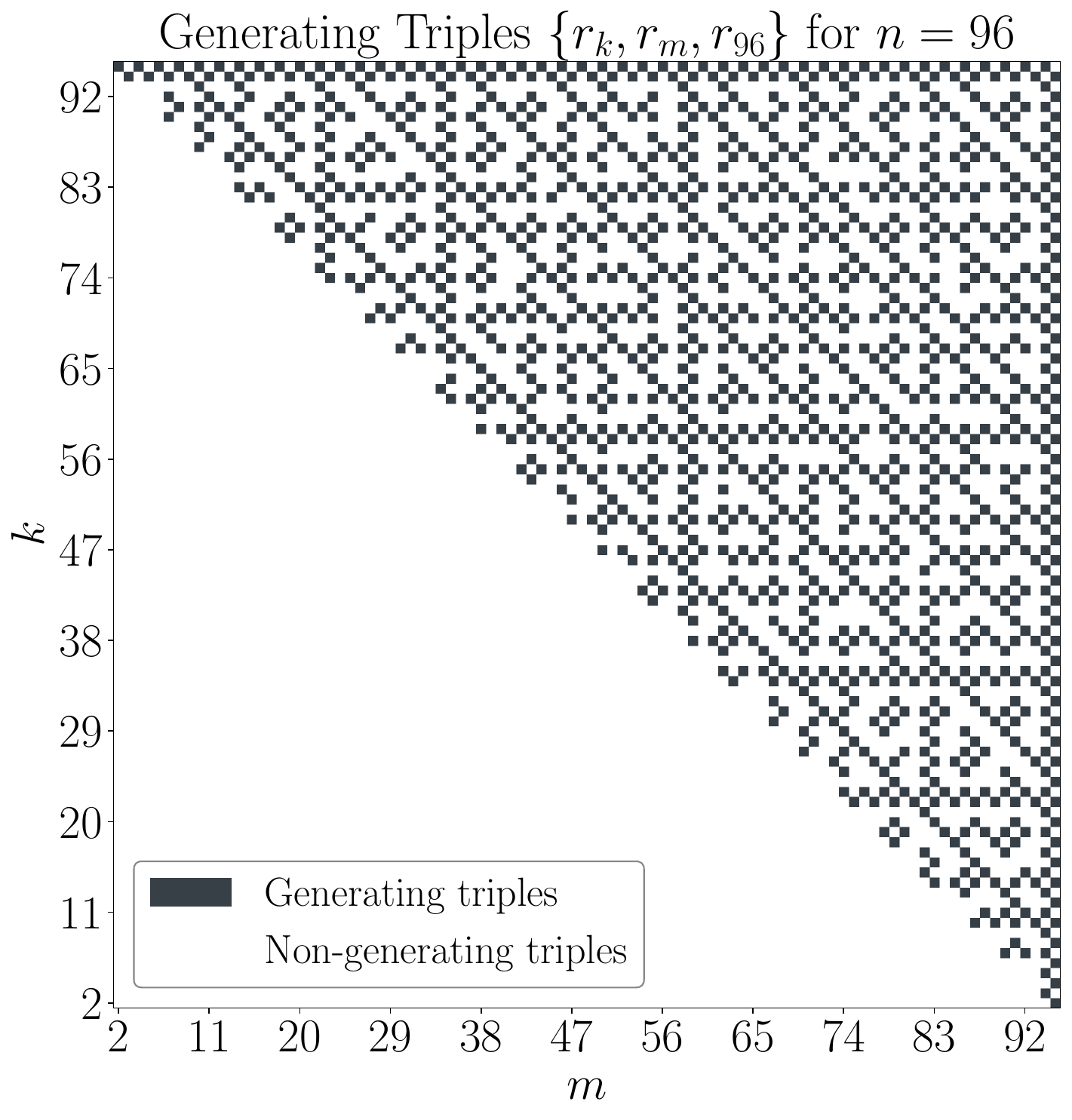}
        \label{fig:97-grid}
    \end{subfigure}
    \hfill
    \begin{subfigure}[b]{0.45\textwidth}
        \includegraphics[width=\textwidth]{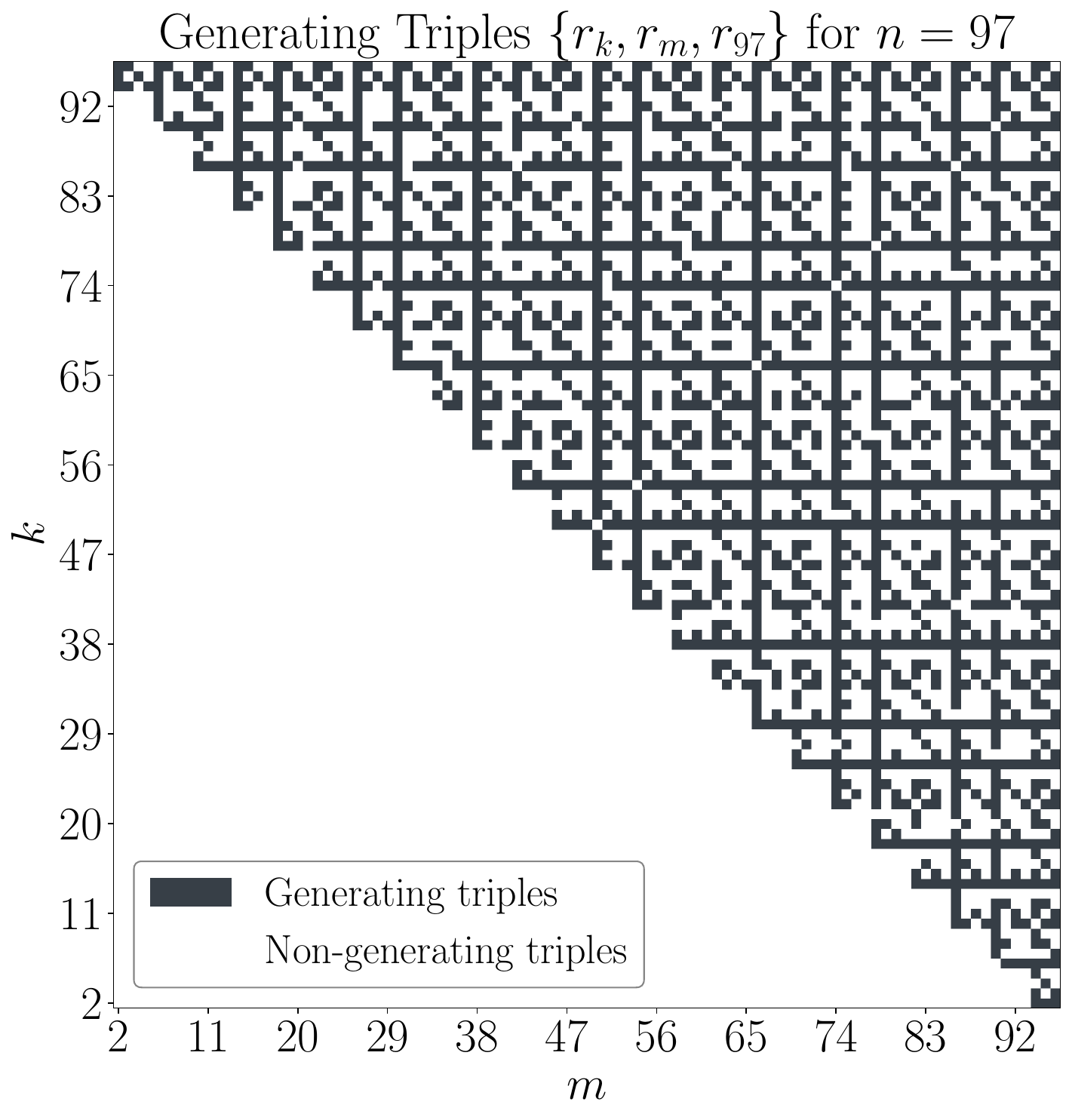}
        \label{fig:98}
    \end{subfigure}
    
    \vspace{0.5cm} 
    
    \begin{subfigure}[b]{0.45\textwidth}
        \includegraphics[width=\textwidth]{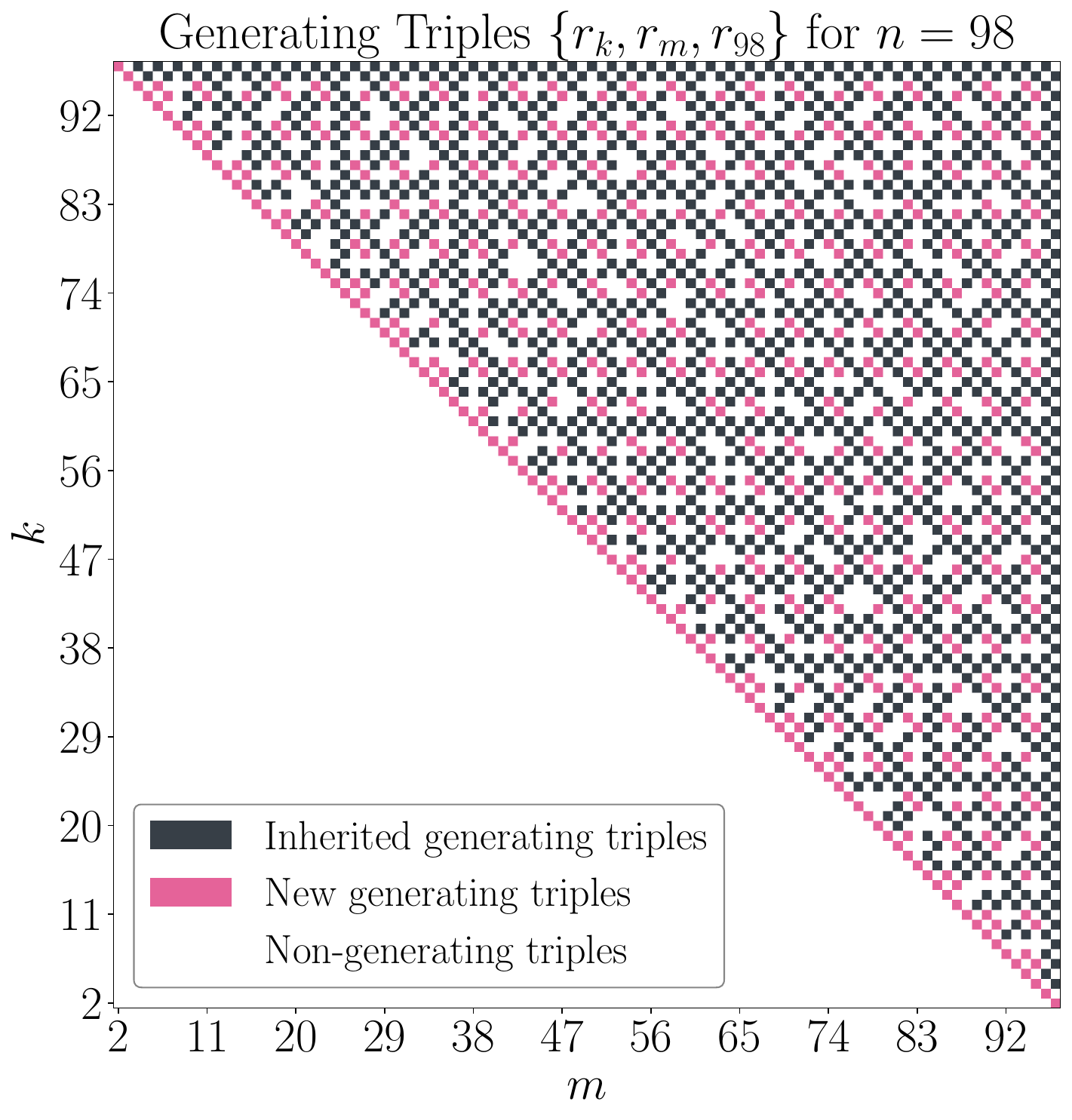}
        \label{fig:99}
    \end{subfigure}
    \hfill
    \begin{subfigure}[b]{0.45\textwidth}
        \includegraphics[width=\textwidth]{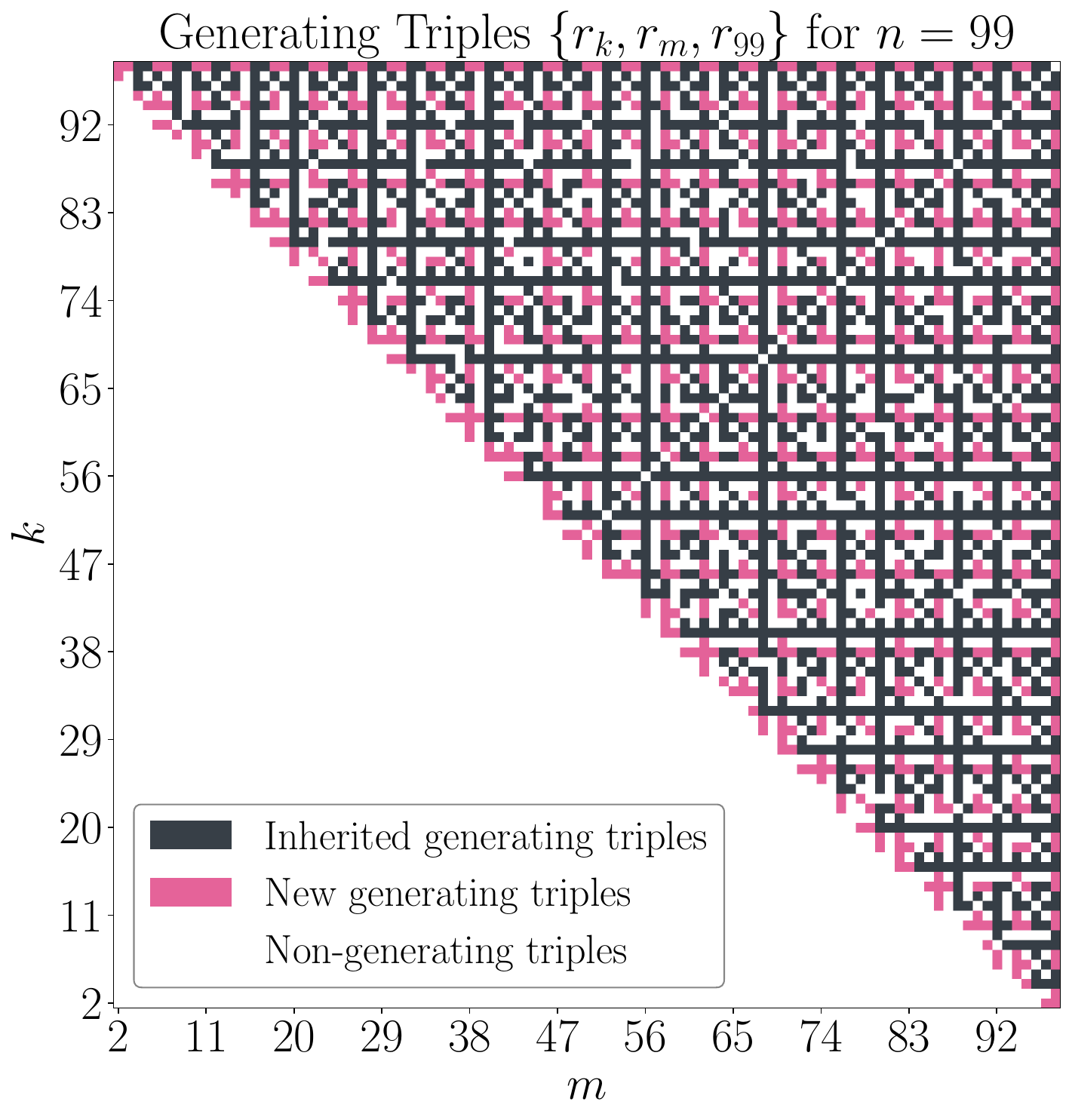}
        \label{fig:100}
    \end{subfigure}
    
    \caption{Visualization of Conjecture~\ref{conj:mod2} to show on how the $n$ case is ``embedded'' to the $n+2$ case. In the first pictures in the left-most column, we depict all generating triples for $n=96$. All of these triples give rise to generating triples for $n=98$ by adding $2$ to each index (depicted in pink). The figures in the right-most column are built in a similar fashion.}
    \label{fig:mod2}
\end{figure}

These conjectures arise naturally by visualizing the set of triples $\{r_n, r_m, r_k\}$ that generate the whole group for a fixed value of $n$. Indeed, in Figure~\ref{fig:mod4} and Figure~\ref{fig:mod2}, we plot the pairs $(m,k)$ that generate $\Sym_n$. Then, to visualize Conjectures~\ref{conj:mod4} and \ref{conj:mod2} in terms of these graphs notice that that triples seem to be preserved by adding the same number to each of the indices. Moreover, Conjecture~\ref{conj:km_sum} can be seen as the white triangle on the plots below the diagonal. 

Both from Figures~\ref{fig:mod4} and \ref{fig:mod2}, one can see that for each ``row'' it is either ``almost full'' or contains only half of the black dots. This observation along with other computational results gives us the following conjecture.  

\begin{conjecture}\label{con:F_k(n)}
    Let 
    $$
    F_k(n) = \{m \,:\, k< m < n, \, \langle r_n, r_m, r_k \rangle = \Sym_n\},
    $$
    then 
    $$
    \max_{n > k+1} |F_k(n)| = \begin{cases}
        k + 1, & \text{if $k$ is even}; \\
        \frac{k + 1}{2}, & \text{if $k$ is odd}. \\
    \end{cases}
    $$
\end{conjecture}

\vspace{5mm}

Taken together, these conjectures suggest that the complete classification may have a periodic and geometric description extending
well beyond the four families settled in this paper.

For the interested reader, the code use to conduct the computations is available at the following link: \url{https://github.com/Luka5s5/prefix_reversal_triples}.

\subsection{Hamiltonicity, diameter, girth}\label{sec:graph}

One of the open problems on Cayley graphs is the conjecture stating that every connected Cayley graph over a finite group has a Hamiltonian cycle. The conjecture was inspired by an open problem about Hamiltonicity of vertex-transitive graphs~\cite{L70}. Hamiltonicity of all $1939864$ cubic Cayley graphs on up to $5000$ vertices is confirmed (see~\url{https://graphsym.net/}).  Computational experiments showed that in the smallest case $n=7$ that has not previously been verified (as $7!>5000$) there are eight possible generating triples $(n,m,k)$ as follows: $(7,4,2)$, $(7,5,2)$, $(7,6,2)$, $(7,4,3)$, $(7,6,3)$, $(7,5,4)$, $(7,6,4)$, and $(7,6,5)$, and the corresponding cubic pancake graphs are Hamiltonian. Therefore, the following statement holds. 

\begin{statement}
Any cubic pancake graph over the symmetric group of degree $n=4,5,6,7$ is Hamiltonian.  
\end{statement}

Moreover, diameters, girths and the shortest cycles were found for all cubic pancake graphs with $n=7$, and the results are summarized in Table~\ref{tab:n=7}. Let us note that Table~\ref{tab:n=7} shows that for $n=7$ the minimal diameter among all the corresponding cubic pancake graphs are associated with the graphs whose generating triples are $(7,6,4)$ and $(7,6,5)$. Meanwhile, girths of these graphs are different; indeed, the graph generated by $(7,6,4)$ has the maximal girth and the graph generated by $(7,6,5)$ has the minimal girth among them.

Additionally, for the Big-3 flips $\{r_n,r_{n-1},r_{n-2}\}$ it was conjectured in~\cite{SW16} that for any $n\geqslant 4$ the corresponding Big-3 pancake graph has cyclic Gray codes, which means it has a Hamiltonian cycle. Due to our computational experiments the following statement holds.

\begin{statement}
The Big-3 pancake graph over the symmetric group of degree $n=4,5,6,7,8$ is Hamiltonian.  
\end{statement}

Totally, for $n=8$ there are five possible triples generating the whole group, namely, $(8,7,2)$,  $(8,6,3)$, $(8,7,4)$, $(8,6,5)$, and $(8,7,6)$. The data on diameters, girths and the shortest cycles of the corresponding graphs are given in Table~\ref{tab:n=8}. 

\begin{table}[h]
\centering
\begin{tabular}{|c|c|c|c|}
\hline
Graph & diameter & girth & shortest cycles\\
\hline
(7,4,2) & 20 & 8 & $(r_4r_2)^4$, $(r_7r_2)^4$,  $(r_7r_4r_7r_2)^2$ \\ 
\hline
(7,5,2) & 20 & 8 & $(r_5r_2)^4$, $(r_7r_2)^4$,  $(r_7r_5r_7r_2)^2$ \\ 
\hline
(7,6,2) & 18 & 8 & $(r_6r_2)^4$, $(r_7r_2)^4$, \\ 
\hline
(7,4,3) & 20 & 8 & $(r_4r_3)^4$, $(r_7r_4r_7r_3)^2$ \\ 
\hline
(7,6,3) & 16 & 8 & $(r_6r_3)^4$ \\ 
\hline
(7,5,4) & 17 & 10 & $(r_5r_4)^5$, $(r_7r_5r_7r_4r_5r_7r_4r_7r_5r_4)$ \\
\hline
(7,6,4) & 15 & 12 & 
\begin{tabular}{c}
$(r_6r_4)^6$, $(r_7r_6r_4)^4$, \\ $(r_7r_6r_7r_4r_7r_4)^2$,  $(r_7r_6r_7r_6r_4r_6)^2$
\end{tabular} \\ \hline
(7,6,5) & 15 & 8 & $(r_7r_6r_5r_6)^2$ \\
\hline
\end{tabular}
\caption{Diameters, girths and shortest cycles for all cubic pancake graphs when $n=7$. Cycles are presented by their canonical and contracted forms.}
\label{tab:n=7}
\end{table}

\begin{table}[h]
\centering
\begin{tabular}{|c|c|c|c|}
\hline
Graph & diameter & girth & shortest cycles\\
\hline
(8,7,2) & 24 & 8 & $(r_7r_2)^4$ \\ \hline
(8,6,3) & 24 & 8 & $(r_6r_3)^4$ \\ \hline
(8,7,4) & 27 & 8 & $(r_8r_4)^4$, $(r_8r_4r_7r_4)^2$ \\ \hline
(8,6,5) & 20 & 12 & $(r_6r_5)^6$, $(r_8r_6r_8r_5r_8r_5)^2$ \\  \hline
(8,7,6) & 16 & 8 & $(r_8r_7r_6r_7)^2$ \\
\hline
\end{tabular}
\caption{Diameters, girths and shortest cycles for all cubic pancake graphs with $n=8$, where cycles are presented by their contracted forms.} 
\label{tab:n=8}
\end{table}

To represent the shortest cycles in Table~\ref{tab:n=7} and Table~\ref{tab:n=8}, we use their {\it canonical forms} $C_{l}=r_{i_0} \ldots r_{i_{l-1}}$ with taking the lexicographically maximal sequence of indices $i_0\ldots i_{\ell-1}$, where $2\leqslant i_j \leqslant n$, and $i_j \neq i_{j+1}$ for any $0\leqslant j \leqslant \ell-1$. {\it Contracted forms} are also used in the tables such that $(r_a r_b)^k$ corresponds to the cycle $C_\ell=r_a r_b\ldots r_a r_b$, where $\ell=2\,k, \ a\neq b$, and $r_a r_b$ appears exactly $k$ times. A similar notation is used for cycles of forms $(r_a\ldots r_c)^2$ and $(r_ar_br_c)^4$, where $(r_a\ldots r_c)^2=(r_a\ldots r_cr_a\ldots r_c)$ and $(r_ar_br_c)^4=(r_ar_br_cr_ar_br_cr_ar_br_cr_ar_br_c)$. For example, $(r_7r_6r_4)^4$ corresponds to a cycle of length $12$ whose canonical form is $C_{12}=(r_7r_6r_4r_7r_6r_4r_7r_6r_4r_7r_6r_4$). All the forms presented in the tables were obtained by Lemma~1 from~\cite{KM16}, Theorem~1.3 from~\cite{KM14}, and computational experiments.

In Table~\ref{tab:minimal_properties}, all generating sets are presented for the cubic pancake graphs that have minimal diameters and girths among all possible cubic pancake graphs with a given $n$, where $4\leqslant n \leqslant 10$. As one can see from the table, for any $n\geqslant 6$, the minimal girth is $8$. These computational results are correlated with theoretical results from~\cite{KS22}. All generating sets of the cubic pancake graphs with maximal diameters and girths are given in Table~\ref{tab:maximal_properties}.

\begin{table}[htbp]
\centering
\small
\setlength{\tabcolsep}{4pt}
\renewcommand{\arraystretch}{1.2}
\begin{tabular}{|c| c| >{\raggedright\arraybackslash}p{2.4cm}| c| >{\raggedright\arraybackslash}p{3.8cm}|}
\hline
$n$ & min diameter & min diameter triples & min girth & min girth triples \\
\hline
4 & 4 & (4,3,2) & 6 & (4,3,2) \newline \\ \hline
5 & 7 & (5,3,2), (5,4,2) & 6 & (5,3,2) \newline \\ \hline
6 & 12 & (6,5,2), (6,4,3), (6,5,4) & 8 & (6,5,2), (6,4,3), (6,5,4) \\ \hline
7 & 15 & (7,6,4), (7,6,5) & 8 & (7,4,2), (7,5,2), (7,6,2), (7,4,3), (7,6,3), (7,6,5) \\ \hline
8 & 20 & (8,6,5) & 8 & (8,7,2), (8,6,3), (8,7,4), (8,7,6) \\ \hline
9 & 23 & (9,6,4), (9,7,6) & 8 & (9,7,2), (9,8,2), (9,8,3), (9,6,4), (9,7,4), (9,8,7) \\ \hline
10 & 27 & (10,9,6), (10,8,7) & 8 & (10,9,2), (10,8,3), (10,9,4), (10,6,5), (10,8,5), (10,9,8) \\
\hline
\end{tabular}
\caption{Minimal diameters and girths among all possible cubic pancake graphs generated by prefix-reversal triples $(n,m,k)$ for a given $n$, where $n \leqslant 10$ and $2\leqslant k < m < n$.}
\label{tab:minimal_properties}
\end{table}

\begin{table}[htbp]
\centering
\small
\setlength{\tabcolsep}{4pt}
\renewcommand{\arraystretch}{1.2}
\begin{tabular}{|c| c| >{\raggedright\arraybackslash}p{2.4cm}| c| >{\raggedright\arraybackslash}p{3.8cm}|}
\hline
$n$ & max diameter & max diameter triples & max girth & max girth triples \\
\hline
4 & 4 & (4,3,2) & 6 & (4,3,2) \newline \\ \hline
5 & 8 & (5,4,3) & 8 & (5,4,2), (5,4,3)\newline \\ \hline
6 & 12 & (6,5,2), (6,4,3), (6,5,4) & 8 & (6,5,2), (6,4,3), (6,5,4) \\ \hline
7 & 20 & (7,4,2), (7,5,2), (7,4,3) & 12 & (7,6,4) \\ \hline
8 & 27 & (8,7,4) & 12 & (8,6,5) \newline \\ \hline
9 & 42 & (9,7,2) & 12 & (9,6,5), (9,7,6), (9,8,6) \newline \\  \hline
10 & 39 & (10,9,2), (10,8,3) & 16 & (10,8,7) \newline \\
\hline
\end{tabular}
\caption{Maximal diameters and girths among all possible cubic pancake graphs generated by prefix-reversal triples $(n,m,k)$ for a given $n$, where $n \leqslant 10$ and $2\leqslant k < m < n$.}
\label{tab:maximal_properties}
\end{table}

\acknowledgements
\label{sec:ack}
All authors contributed equally and are listed alphabetically by their names in English. The authors thank the anonymous reviewers for their valuable comments that improved the presentation of this paper.

Elena V.~Konstantinova thanks the N.N. Krasovskii Institute of Mathematics and Mechanics and Ural Mathematical Center that hosted her visit to Yekaterinburg in $2023$ when this project was initiated. She was supported by the Mathematical Center in Akademgorodok under the agreement No.~075-15-2025-348 with the Ministry of Science and Higher Education of the Russian Federation. Mikhail Golubyatnikov and  Natalia Maslova worked within the framework of the collaboration between the Ural Mathematical Center and the Mathematical Center in Akademgorodok. Sa\'ul A.~Blanco thanks Elena Konstantinova for her kind hospitality during the 2025 G2A2 Conference in Novosibirsk where part of this work was discussed. 


\newpage

\nocite{*}
\bibliographystyle{abbrvnat}
\bibliography{reference.bib}
\label{sec:biblio}

\end{document}